\documentclass[12pt,twoside]{amsart}
\usepackage{amssymb}
\usepackage[ngerman, english]{babel}
\usepackage{amscd}
\usepackage{hyperref}
\usepackage{color}

\usepackage[normalem]{ulem} 
\newcommand\redout{\bgroup\markoverwith
{\textcolor{red}{\rule[.5ex]{2pt}{0.4pt}}}\ULon}

\title[Vanishing theorems of Kodaira type]
{Vanishing theorems of Kodaira type for Witt canonical sheaves} 
\author{Hiromu Tanaka} 
\address{Graduate School of Mathematical Sciences, 
The University of Tokyo, 
3-8-1 Komaba, Meguro-ku, Tokyo 153-8914, JAPAN} 
\email{tanaka@ms.u-tokyo.ac.jp}

\newcommand{\Ker}[0]{{\operatorname{Ker}}}

\newcommand{\Spec}[0]{{\operatorname{Spec}}}

\newcommand{\Supp}[0]{{\operatorname{Supp}}}

\newtheorem{thm}{Theorem}[section]
\newtheorem{lem}[thm]{Lemma}
\newtheorem{cor}[thm]{Corollary}
\newtheorem{prop}[thm]{Proposition}

\newtheorem{step}{Step}

\theoremstyle{definition}

\newtheorem{dfn}[thm]{Definition}

\newtheorem{rem}[thm]{Remark}

\newtheorem{nothing}[thm]{}

\makeatletter
  
  \@addtoreset{equation}{thm}
  \makeatother

\newcommand{\MO}{\mathcal{O}}

\newcommand{\F}{\mathbb{F}}

\newcommand{\R}{\mathbb{R}}
\newcommand{\Q}{\mathbb{Q}}
\newcommand{\Z}{\mathbb{Z}}

\newcommand{\m}{\mathfrak{m}}

\begin{document}

\maketitle

\begin{abstract}
Given a smooth projective variety over a perfect field of positive  characteristic, we prove that the higher cohomologies vanish for the tensor  product of the Witt canonical sheaf and the Teichm\"{u}ller lift of an ample  invertible sheaf. We also give a generalisation of this vanishing theorem to one of Kawamata-Viehweg type. 
\end{abstract}

\tableofcontents

\section{Introduction}

One of deep results of complex algebraic geometry 
is the degeneration of Hodge to de Rham spectral sequence. 
Indeed, it is related to several important results such as the Hodge decomposition and 
the Kodaira vanishing theorem. 
Unfortunately Kodaira vanishing is known to fail in positive characteristic \cite{Ray78}. 
On the other hand, 
one of positive-characteristic counterparts of the Hodge decomposition 
is the slope decomposition of the Crystalline cohomologies by 
the de Rham--Witt complex \cite{Ill79}. 
Thus it is natural to seek a positive-characteristic analogue of the Kodaira vanishing theorem in terms of de Rham--Witt complex. 
One of the purposes of this paper is to propose candidates as follows. 

\begin{thm}[cf. Theorem \ref{t-Kod-van}, Theorem \ref{t-KVV-minus}]\label{intro-main1}
Let $k$ be a perfect field of characteristic $p>0$ and 
let $X$ be an $N$-dimensional smooth projective connected scheme over $k$. 
Let $A$ be an ample invertible sheaf on $X$ and let $\underline{A}$ be 
the Teichm\"{u}ller lift of $A$ (cf. Definition \ref{d-teich-lift}). 
Then 
\begin{enumerate}
\item $H^i(X, W\Omega_X^N \otimes_{W\MO_X} \underline{A})=0$ 
for any positive integer $i$, and 
\item 
$H^j(X, \underline{A}^{-1}) \otimes_{\Z} \Q=0$ for any integer $j$ satisfying $j<N$. 
\end{enumerate}
\end{thm}

\begin{rem}\label{intro-r-Q-Kodaira}
We actually need to take the tensor product with $\Q$ in the statement (2) of Theorem \ref{intro-main1}, as $H^j(X, \underline{A}^{-1})$ does not vanish in general 
(cf. Proposition \ref{p-cex-torsion}).
\end{rem}

Since we have established an analogue of the Kodaira vanishing theorem, 
it is tempting to find analogous results related to the Kodaira vanishing theorem. 
In this direction, we shall give a generalisation of Theorem \ref{intro-main1} 
that can be considered as a vanishing theorem of Kawamata--Viehweg type.

\begin{thm}[cf. Corollary \ref{c-KVV}, Theorem \ref{t-KVV-minus}]\label{intro-main2}
Let $k$ be a perfect field of characteristic $p>0$. 
Let $f:X \to Y$ be a projective $k$-morphism 
from an $N$-dimensional connected scheme $X$ smooth over $k$ 
to a scheme $Y$ of finite type over $k$. 
Let $A$ be an $f$-ample $\Q$-divisor on $X$ 
such that $\Supp\,\{A\}$ is simple normal crossing. 
Let $W\MO_X(-A)$ be the Witt divisorial sheaf associated to $-A$ 
(cf. Definition \ref{d-div-def}).  
Then 
\begin{enumerate}
\item 
$R^if_*(\mathcal Hom_{W\MO_X}(W\MO_X(-A), W\Omega_X^N)) \otimes_{\Z} \Q=0$ 
for any positive integer $i$, and  
\item 
if the induced morphism $Y \to \Spec\,k$ is an isomorphism, then 
$$H^j(X, W\MO_X(-A)) \otimes_{\Z} \Q=0$$ 
for any integer $j$ satisfying $j<N$. 
\end{enumerate}
\end{thm}

\begin{rem}
One of the forms of the Kawamata--Viehweg vanishing theorem in characteristic zero 
is the following vanishing 
\[
R^if_*(\mathcal Hom_{\MO_X}(\MO_X(-A), \Omega_X^N)) = 0
\]
under the same assumption as in Theorem \ref{intro-main2} except for the characteristic of $k$. 
Indeed, this is the relative version of \cite[Corollary 1-2-2]{KMM87}. 
Therefore, Theorem \ref{intro-main2} can be considered as a Witt analogue of the Kawamata--Viehweg vanishing theorem. 
\end{rem}

\begin{rem}
If $A$ is a $\Z_{(p)}$-divisor, then 
(1) of Theorem \ref{intro-main2} holds even without $(-) \otimes_{\Z} \Q$ (Theorem \ref{t-KVV}). 
The author does not know whether we may drop $(-) \otimes_{\Z} \Q$ in general. 
\end{rem}

\begin{rem}
We now summarise open problems related to Theorem \ref{intro-main1} and Theorem \ref{intro-main2}. 
Let $k$ be a perfect field of characteristic $p>0$ and 
let $X$ be an $N$-dimensional smooth projective connected scheme over $k$. 
\begin{enumerate}
\item 
If $A$ is a nef and big invertible sheaf on $X$, then 
does it hold that $H^i(X, W\Omega_X^N \otimes_{W\MO_X} \underline{A}) \otimes_{\Z} \Q=0$ for $i > 0$? 
\item 
If $A$ is an ample invertible sheaf on $X$, then 
does it hold that $H^i(X, W\Omega_X^j \otimes_{W\MO_X} \underline{A})=0$ for $i+j > \dim X$? 
\end{enumerate}
The original motivation of the author was to prove (1), as its relative version 
deduces the Grauert--Riemenschneider vanishing for Witt canonical sheaves. 
Contrary to the situation in characteristic zero, Theorem \ref{intro-main2} does not seem to be enough to conclude (1). 
\end{rem}


It is worth mentioning that 
the cohomologies of the Teichm\"{u}ller lifts of invertible sheaves 
can be infinite dimensional over $W(k) \otimes_{\Z} \Q$. 
More specifically, we shall prove the following statement.

\begin{thm}[cf. Theorem \ref{t-ample-dim},  Theorem \ref{t-highest-infinite}]\label{intro-main3}
Let $k$ be a perfect field of characteristic $p>0$ and set $Q:=W(k) \otimes_{\Z} \Q$. 
Let $X$ be a projective scheme over $k$ such that $\dim X>0$. 
Let $A$ be an ample invertible sheaf on $X$ and let $\underline{A}$ be 
the Teichm\"{u}ller lift of $A$ (cf. Definition \ref{d-teich-lift}). 
Then the following hold. 
\begin{enumerate}
\item 
$H^i(X, \underline{A}) \otimes_{\Z} \Q=0$ for any $i>0$. 
\item 
$\dim_Q (H^0(X, \underline{A}) \otimes_{\Z} \Q)=\infty$. 
\item 
If $X$ is smooth over $k$, then $\dim_Q (H^{\dim X}(X, \underline{A}^{-1}) \otimes_{\Z} \Q)=\infty$. 
\end{enumerate}
\end{thm}

\subsection{Description of proofs}

\subsubsection{The proof of Theorem \ref{intro-main1}}\label{sss-main1}

Let us overview some of the ideas of the proof of Theorem \ref{intro-main1}(1). 
The argument consists of three steps. 
The first step is to twist the finite levels $W_n(\bullet)$ by Frobenius 
in order to achieve the vanishing, i.e. we can find a positive integer $t$ such that 
\[
H^i(X, W_n\Omega_X^N \otimes_{W\MO_X} (F_X^t)^*\underline{A}) 
\simeq 
H^i(X, W_n\Omega_X^N \otimes_{W\MO_X} \underline{A}^{p^t})=0
\]
for any $i>0$ and $n>0$, where $F_X:X \to X$ is the absolute Frobenius morphism. 
We can find such a positive integer $t$ by the exact sequence (cf. the proof of Theorem \ref{t-Kod-van}(2)): 
\[
0 \to \Omega_X^N \to W_{n+1}\Omega_X^N \to (F_X)_*(W_n\Omega_X^N) \to 0. 
\]
The second step is to go to the infinite level $W(\bullet)$ from the finite levels $W_n(\bullet)$, 
i.e. we take the projective limit: 
\[
H^i(X, W\Omega_X^N \otimes_{W\MO_X} (F_X^t)^*\underline{A}) \simeq 
\varprojlim_n H^i(X, W_n\Omega_X^N \otimes_{W\MO_X} (F_X^t)^*\underline{A})=0, 
\]
where the above isomorphism holds by the Mittag--Leffler condition. 
The final step is to untwist the infinite level $W(\bullet)$ by Frobenius, i.e. 
\begin{eqnarray*}
H^i(X, W\Omega_X^N \otimes_{W\MO_X} \underline{A}) 
&\simeq & H^i(X, (F^t_X)_*(W\Omega_X^N) \otimes_{W\MO_X} \underline{A})\\
&\simeq & H^i(X, (F^t_X)_*(W\Omega_X^N \otimes_{W\MO_X} (F^t_X)^*\underline{A}))\\
&\simeq & H^i(X, W\Omega_X^N \otimes_{W\MO_X} (F^t_X)^*\underline{A})\\
&=&0. 
\end{eqnarray*}
We can apply such an untwisting argument by using the $W\MO_X$-linear isomorphism 
$W\Omega_X^N \simeq (F^t_X)_*(W\Omega_X^N)$ (Theorem \ref{t-descent}). 

\subsubsection{The relative version of Theorem \ref{intro-main1}}

In Theorem \ref{t-Kod-van}, we shall establish the relative version of Theorem \ref{intro-main1}. 
The proof is identical to the above argument in Subsubsection \ref{sss-main1} 
except for the second step, i.e. the step taking the projective limit. 
This is because we can not use the Mittag--Leffler criterion any more. 
For the details on this argument, see the proof of Theorem \ref{t-Kod-van}(3). 

\subsubsection{Sketch of the proof of Theorem \ref{intro-main2}}

We now overview the outline of the proof of Theorem \ref{intro-main2}(1). 
The main idea is to apply Kawamata's covering trick, 
as in Kawamata's original proof of the Kawamata--Viehweg vanishing theorem. 
To this end, we first take a base change $f:X':=X \times_k k' \to X$ 
for a suitable finite extension $k \subset k'$. 
Applying Kawamata's covering trick \cite[Theorem 1-1-1]{KMM87}, 
we obtain a finite Galois cover $g:X'' \to X'$ 
such that $[K(X''):K(X')]$ is not divisible by $p$ and 
the pullback of $A$ to $X''$ is a $\Z$-divisor. 
Let $\varphi:Z' \to Z$ be one of $f$ and $g$. 
For $V \in \{Z, Z'\}$, we denote by $A_V$ the pullback of $A$ to $V$. 
Roughly speaking, the problem is reduced to prove the following two assertions: 
\begin{enumerate}
\item[(I)] $W\MO_{Z}(-A_Z) \to \varphi_*W\MO_{Z'}(-A_{Z'})$ splits. 
\item[(II)] The induced homomorphism 
\[
\mathcal Hom_{W_n\MO_X} (\varphi_*W_n\MO_{Z'}(-A_{Z'}), W_n\Omega_Z^N) \to 
\mathcal Hom_{W\MO_X}(\varphi_*W\MO_{Z'}(-A_{Z'}), W_n\Omega_Z^N)
\]
is an isomorphism. 
\end{enumerate}

The proof of (I) is rather standard. 
More specifically, we shall check that 
several standard arguments on $\MO_X(D)$ can be applied for $W\MO_X(D)$. 
For more details, see Subsection \ref{ss4-split}. 

Subsection \ref{ss-WO-WO} is devoted to prove (II), 
which is the most technical part in the proof of Theorem \ref{intro-main2}. 
We introduce notion of weakly $\ell$-cyclic morphisms  
in order to treat the two cases $\varphi=f$ and $\varphi=g$ simultaneously (Definition \ref{d-weak-l-cyclic}). 
For the technical details, see the proofs of Lemma \ref{l-WO-vs-WO-Y} and Proposition \ref{p-WO-vs-WO-Y}.

\subsection{Witt divisorial sheaves}

Another novelty of this paper is the notion of Witt divisorial sheaves 
$W\MO_X(D)$ and $W_n\MO_X(D)$ for an $\R$-divisor $D$. 
Section \ref{s-Witt-div-sh} is devoted to establishing their fundamental properties. 
If $D$ is a Cartier divisor, then 
the Witt divisorial sheaf $W\MO_X(D)$ coincides with the Teichm\"{u}ller lift of 
the invertible sheaf $\MO_X(D)$ (Proposition \ref{p-div-Teich}). 
The author introduces this notion to find a Witt analogue of 
the Kawamata--Viehweg vanishing theorem (Theorem \ref{intro-main2}). 
It is remarkable that $W_n\MO_X(D)$, for certain special cases, has already appeared in the study of Swan conductor by Brylinski--Kato (cf. Remark \ref{r-BK}). 
In this sense, the Witt divisorial sheaves naturally appear in different situations. 
The author hopes that $W\MO_X(D)$ will produce further applications in future.


\medskip

\textbf{Acknowledgements:} 
The author would like to thank 
Piotr Achinger, 
Yoshinori Gongyo, 
Luc Illusie, 
Yujiro Kawamata, 
Yukiyoshi Nakkajima, 
and Shuji Saito 
for useful comments and pointing out mistakes. 
The author also would like to thank the referee for 
reading the manuscript carefully and suggesting several improvements. 
The author was funded by EPSRC 
and the Grant-in-Aid for Scientific Research (KAKENHI No. 18K13386).

\section{Preliminaries}

\subsection{Notation}\label{ss-notation}

In this subsection, we summarise notation used in this paper. 

\begin{enumerate}
\item 
Throughout the paper, $p$ denotes a prime number. 
\item We shall freely use the notation and terminology in \cite{Har77} 
and \cite{Ill79}. 
\item Let $N$ be a non-negative integer. 
A noetherian scheme $X$ is {\em pure} $N$-{\em dimensional} 
if $\dim Y=N$ for any irreducible component $Y$ of $X$. 
A scheme $X$ is {\em excellent} if 
$X$ is a noetherian scheme such that $\MO_{X, x}$ is an excellent ring for any $x \in X$.  
\item 
For an $\F_p$-scheme $X$, 
we denote by $F_X:X \to X$ the {\em absolute Frobenius morphism}, 
i.e. the homeomorphic affine morphism such that 
the induced ring homomorphism can be written by 
$$\Gamma(U, \MO_X) \to \Gamma(U, \MO_X), \quad a \mapsto a^p$$
for any open subset $U$ of $X$. 
An $\F_p$-scheme $X$ is $F$-{\em finite} if $F_X:X \to X$ is a finite morphism. 
\item For an integral scheme $X$, 
we define the {\em function field} $K(X)$ of $X$ 
as the stalk $\MO_{X, \xi}$ at the generic point $\xi$ of $X$. 
\item For a field $k$, 
$X$ is a {\em variety over} $k$ or a $k$-{\em variety} if 
$X$ is an integral scheme that is separated and of finite type over $k$. 
$X$ is a {\em curve} over $k$ or a $k$-{\em curve} 
(resp. a {\em surface} over $k$ or a $k$-{\em surface}) 
if $X$ is a $k$-variety of dimension one (resp. two). 
\item 
Let $X$ be a regular noetherian scheme. 
A closed subset $Z$ of $X$ is {\em simple normal crossing} if 
for the irreducible decomposition $Z=\bigcup_{i \in I} Z_i$ of $Z$ and any subset $J$ of $I$, the scheme-theoretic intersection $\bigcap_{j \in J}Z_j$ is a regular scheme, 
where each $Z_i$ is equipped with the reduced scheme structure. 
\item  
A morphism $f:X \to S$ of noetherian schemes 
is {\em projective} if there exists a closed immersion 
$X \hookrightarrow \mathbb P^n_S$ over $S$ for some non-negative integer $n$. 
We adopt the definition of projective morphisms by \cite{Har77} 
(cf. \cite[Section 5.5.1]{FGAex}). 
\item\label{ss-notation-9} 
For an integral normal noetherian scheme $X$ and a subring $\mathbb K$ of $\mathbb R$, 
an $\mathbb K$-{\em divisor} $D$ on $X$ is an $\mathbb K$-linear combination 
$\sum_{i \in I} a_i D_i$, 
i.e. $I$ is a finite set, $a_i \in \mathbb K$, and $D_i$ is a prime divisor 
for any $i \in I$.  
An $\mathbb R$-divisor $D$ is $\mathbb K$-{\em Cartier} 
if there is an equation $D=\sum_{j=1}^s \alpha_j E_j$ 
for some $\alpha_1, \cdots, \alpha_s \in \mathbb K$ and Cartier divisors $E_1, \cdots, E_s$. 
Note that a $\mathbb K$-Cartier $\R$-divisor is a $\mathbb K$-divisor. 
Given a dominant morphism $f:Y \to X$ of integral normal noetherian schemes 
and a $\mathbb K$-Cartier $\mathbb K$-divisor $D$ with the equation $D=\sum_{j=1}^s \alpha_j E_j$ as above, we define the {\em pullback} $f^*D$ of $D$ 
as $\sum_{j=1}^s \alpha_j f^*E_j$. 
We can check that this does not depends on 
the choice of expression $D=\sum_{j=1}^s \alpha_j E_j$. 
\item 
For $r \in \R$, we define $\llcorner r\lrcorner$ (resp. $\ulcorner r\urcorner$) as the integer such that 
$\llcorner r \lrcorner \leq r <\llcorner r \lrcorner+1$ 
(resp. $\ulcorner r\urcorner-1< r\leq \ulcorner r\urcorner$). 
For an $\R$-divisor $D$ on an integral normal noetherian scheme $X$ 
and its irreducible decomposition $D=\sum_{i \in I} r_i D_i$, 
we set $\llcorner D\lrcorner:=\sum_{i \in I} \llcorner r_i \lrcorner D_i$, 
$\ulcorner D\urcorner:=\sum_{i \in I} \ulcorner r_i \urcorner D_i$, 
and $\{D\}:=D-\llcorner D\lrcorner$. 
\item 
Let $f:X \to S$ be a projective morphism of noetherian schemes. 
An invertible sheaf $L$ on $X$ is $f$-{\em ample} 
if there exists a non-negative integer $n$, 
a closed immersion $j:X \hookrightarrow \mathbb P^n_S$ over $S$, 
and a positive integer $m$ 
such that $L^{\otimes m} \simeq j^*\MO_{\mathbb P^n_S}(1)$. 
A Cartier divisor $D$ is $f$-{\em ample} if 
the invertible sheaf $\MO_X(D)$ is $f$-ample. 
Assume that $X$ is an integral normal scheme. 
An $\R$-Cartier $\R$-divisor $E$ on $X$ is $f$-{\em ample} 
if we can write $E=\sum_{j \in J} a_j D_j$ for some non-empty finite set $J$, 
$a_j \in \R_{>0}$, 
and $f$-ample Cartier divisors $D_j$. 
An $\R$-Cartier $\R$-divisor $D$ is {\em ample} if $D$ is $g$-ample for some projective morphism $g:X \to T$ to a noetherian affine scheme $T$. 
\item 
For an $\F_p$-scheme $X$, 
we denote the ringed space $(X, W_n\MO_X)$ by $W_nX$. 
It is well known that $W_nX$ is a scheme. 
If $X$ is an $F$-finite noetherian scheme, then $W_nX$ is a noetherian scheme 
\cite[Ch. 0, the discussion after (1.5.3) in page 512]{Ill79}. 
Furthermore, the morphism $W_n(F_X): W_nX \to W_nX$ induced by the absolute Frobenius morphism 
$F_X:X \to X$ is a finite morphism \cite[Proposition 1.5.6(ii)]{Ill79}. By abuse of notation, 
$W_n(F_X)$ is denoted by $F_X$. 
\item 
Let $X$ be an $\F_p$-scheme. 
It is well known that we have the natural surjective closed immersion $X \to W_nX$ 
corresponding to $W_n\MO_X \to \MO_X, (a_0, a_1, ...) \mapsto a_0$. 
Since $X \to W_nX$ is a homeomorphism, we often identify these topological spaces. 
We have the {\em Teichm\"uller lift} $\MO_X \to W_n\MO_X, a \mapsto (a, 0, 0, ...) =: \underline{a}$, 
which is known to be a multiplicative map, i.e. 
the equation $\underline{a} \cdot \underline{b} = \underline{ab}$ 
holds for any open subset $U$ of $X$ and $a, b \in \MO_X(U)$. 
\item\label{ss-notation-Q} 
For an abelian group $H$, we set $H_{\Q}:=H \otimes_{\Z} \Q$. 
For a sheaf $F$ of abelian groups on a topological space $X$, 
$F_{\Q}$ is the presheaf defined by $\Gamma(U, F_{\Q}) := \Gamma(U, F)_{\Q}$ for any open subset $U$ of $X$. 
Assume that $X$ is a noetherian space. 
It is known that $F_{\Q}$ is automatically a sheaf.  
Furthermore, we have an isomorphism $H^i(X, F)_{\Q} \simeq H^i(X, F_{\Q})$. 
Indeed, given a flasque resolution $0 \to F \to I^0 \to I^1 \to \cdots$, 
the induced resolution $0 \to F_{\Q} \to I^0_{\Q} \to I^1_{\Q} \to \cdots$ is flasque, 
i.e. each $I^n_{\Q}$ is a flasque sheaf. 
Therefore, we have that 
\[
H^i(X, F)_{\Q} \simeq H^i( \Gamma(X, I^{\bullet} ))_{\Q} \simeq 
H^i( \Gamma(X, I^{\bullet}_{\Q} )) \simeq H^i(X, F_{\Q}), 
\]
as required. 
\end{enumerate}

\subsection{Iterated Cartier operators}

In this subsection, 
we recall definition of $B_n\Omega_X^r$ and $Z_n\Omega_X^r$. 
Our definition slightly differs from the one of \cite[Ch. 0, Section 2.2]{Ill79} (Remark \ref{r-BZ}). 

Let $k$ be a perfect field of characteristic $p>0$ and 
let $X$ be a smooth scheme over $k$. 
It is easy to check that the Frobenius push-forward of 
the de Rham complex: 
$$0 \to (F_X)_*\MO_X \xrightarrow{d_0} (F_X)_*\Omega^1_X \xrightarrow{d_1} (F_X)_*\Omega^2_X \xrightarrow{d_2} \cdots$$
is a complex of $\MO_X$-module homomorphisms, 
where $\Omega^r_X:=\bigwedge^r \Omega^1_{X/\Z}$. 
Since $k$ is a perfect field, 
it holds that $\bigwedge^r \Omega^1_{X/\Z} \simeq \bigwedge^r \Omega^1_{X/k}$. 
For any $r \in \Z_{\geq 0}$, we set 
$$B_1\Omega^r_X:={\rm Im}\,(d_{r-1}), \quad Z_1\Omega^r_X:=\Ker\,(d_r).$$
In particular, both $B_1\Omega^r_X$ and $Z_1\Omega^r_X$ 
are coherent $\MO_X$-submodules of $(F_X)_*\Omega_X^r$. 
It is well known that there is an exact sequence of $\MO_X$-module homomorphisms: 
$$0 \to B_1\Omega^r_X \to Z_1\Omega^r_X \xrightarrow{C} \Omega_X^r \to 0.$$
Fix $n \in \Z_{>0}$ and assume that we have already defined 
$B_n\Omega^r_X$ and $Z_n\Omega^r_X$ as coherent $\MO_X$-submodules of $(F_X^n)_*\Omega_X^r$. 
We have 
$$0 \to (F^n_X)_*B_1\Omega^r_X \to (F^n_X)_*Z_1\Omega^r_X 
\xrightarrow{C_n} (F^n_X)_*\Omega_X^r \to 0,$$
where $C_n:=(F^n_X)_*C$. We set 
$$B_{n+1}\Omega_X^r:=(C_n)^{-1}(B_n\Omega_X^r), \quad 
Z_{n+1}\Omega_X^r:=(C_n)^{-1}(Z_n\Omega_X^r).$$
In particular, both $B_{n+1}\Omega_X^r$ and $Z_{n+1}\Omega^r_X$ are 
coherent $\MO_X$-submodules of $(F_X^{n+1})_*\Omega^r_X$. 
For convenience, we set 
$$B_0\Omega_X^r:=0, \quad Z_0\Omega_X^r:=\Omega_X^r.$$

\begin{rem}\label{r-BZ}
We define $B_n\Omega_X^r$ and $Z_n\Omega_X^r$ 
in a slightly different way from the one of \cite[Ch. 0, Section 2.2]{Ill79}. 
Two definitions are the same up to the base change 
induced by the absolute Frobenius morphism of the base field 
$F_k^n:\Spec\,k \to \Spec\,k$. 
In other words, we treat the absolute version, whilst \cite{Ill79} considers the relative one. 
\end{rem}

\begin{prop}\label{p-BZ}
Let $k$ be a perfect field of characteristic $p>0$ and 
let $X$ be a smooth scheme over $k$. 
Then the following hold. 
\begin{enumerate}
\item 
For any $n \in \Z_{\geq 0}$ and $r \in \Z_{\geq 0}$, 
the sheaves $B_n\Omega_X^r$ and 
$Z_n\Omega_X^r$ are coherent locally free $\MO_X$-modules. 
\item 
For any $n \in \Z_{\geq 0}$ and $r \in \Z_{\geq 0}$, 
there are exact sequences of $\MO_X$-modules: 
$$0 \to (F^n_X)_*B_1\Omega_X^r \to B_{n+1}\Omega_X^r \to B_n\Omega_X^r \to 0,$$
$$0 \to (F^n_X)_*B_1\Omega_X^r \to Z_{n+1}\Omega_X^r \to Z_n\Omega_X^r \to 0.$$
\item 
$(F_X^n)_*\Omega_X^r/ B_n\Omega_X^r$ and 
$(F_X^n)_*\Omega_X^r/ Z_n\Omega_X^r$ are coherent locally free $\MO_X$-modules. 
\end{enumerate}
\end{prop}

\begin{proof}
The assertion (1) holds by \cite[Ch. 0, Proposition 2.2.8(a)]{Ill79} (cf. Remark \ref{r-BZ}). 
The assertion (2) follows from the construction.

Let us show (3) by induction on $n$. 
We first treat the case when $n=1$. 
By the exact sequence
\[
0 \to Z_1\Omega_X^r \to (F_X)_*\Omega_X^r \xrightarrow{d} B_1\Omega_X^{r+1} \to 0, 
\]
it follows from (1) that $(F_X)_*\Omega_X^r/Z_1\Omega_X^r$ is a coherent locally free $\MO_X$-module. 
By the exact sequence 
\[
0 \to B_1\Omega^r_X \to Z_1\Omega^r_X \xrightarrow{C} \Omega_X^r \to 0, 
\]
also $Z_1\Omega^r_X/B_1\Omega^r_X$ is a coherent locally free $\MO_X$-module. 
Furthermore, we have the following exact sequence 
\[
0 \to \frac{Z_1\Omega^r_X}{B_1\Omega^r_X} \to \frac{(F_X)_*\Omega_X^r}{B_1\Omega^r_X} 
\to \frac{(F_X)_*\Omega_X^r}{Z_1\Omega^r_X} \to 0. 
\]
Since an extension of coherent locally free $\MO_X$-modules is a coherent locally free $\MO_X$-module, 
also $(F_X)_*\Omega_X^r/B_1\Omega^r_X$ is a coherent locally free $\MO_X$-module. 
This completes the proof for the case when $n=1$.

Assume that 
$(F_X^n)_*\Omega_X^r/ B_n\Omega_X^r$ and 
$(F_X^n)_*\Omega_X^r/ Z_n\Omega_X^r$ are coherent locally free $\MO_X$-modules. 
We have the following commutative diagram in which each horizontal sequence is exact: 
\[
\begin{CD}
0 @>>> (F_X^n)_*B_1\Omega^r_X @>>> (F_X^n)_*Z_1\Omega^r_X @>C_n>> (F_X^n)_*\Omega_X^r @>>> 0\\
@. @| @AA{\rm inclusion}A @AA{\rm inclusion}A\\
0 @>>> (F_X^n)_*B_1\Omega^r_X @>>> Z_{n+1}\Omega^r_X @>C_n>> Z_n\Omega_X^r @>>> 0.\\
\end{CD}
\]
Since $(F_X^n)_*\Omega_X^r/ Z_n\Omega_X^r$ is a coherent locally free $\MO_X$-module, 
the snake lemma implies that $(F_X^n)_*Z_1\Omega^r_X / Z_{n+1}\Omega^r_X$ 
is a coherent locally free $\MO_X$-module. 
By the case of $n=1$, 
\[
(F_X^{n+1})_*\Omega_X^r / (F_X^n)_*Z_1\Omega^r_X  = 
(F_X^n)_* ( (F_X)_*\Omega_X^r / Z_1\Omega^r_X)
\] 
is a coherent locally free $\MO_X$-module. 
Hence, an extension $(F_X^{n+1})_*\Omega_X^r / Z_{n+1}\Omega^r_X$ of $(F_X^{n+1})_*\Omega_X^r / (F_X^n)_*Z_1\Omega^r_X$ and $(F_X^n)_*Z_1\Omega^r_X / Z_{n+1}\Omega^r_X$ is a coherent locally free $\MO_X$-module. 
By a similar argument, we can show that $(F_X^{n+1})_*\Omega_X^r / B_{n+1}\Omega^r_X$ is a coherent locally free $\MO_X$-module. 
\end{proof}

\subsection{Witt vectors}

For definition and basic properties of Witt vectors, 
we refer to \cite[Ch. 0, Section 1]{Ill79}, \cite[Ch. II, \S 6]{Ser79}. 
In this subsection, we recall two results on Witt vectors (Lemma \ref{l-sum-def}, Lemma \ref{l-sum-homog}) 
for later usage.

\begin{lem}\label{l-sum-def}
There exist unique polynomials 
$$S_0(x_0, y_0), S_1(x_0, y_0, x_1, y_1), \cdots \in \Z[x_0, y_0, x_1, y_1, \cdots]$$
that satisfy the following properties.   
\begin{enumerate}
\item 
For any $n \in \Z_{\geq 0}$, it holds that 
{\small 
$$w_n(S_0(x_0, y_0), \cdots, S_n(x_0, y_0, \cdots, x_n, y_n))=w_n(x_0, \cdots, x_n)+w_n(y_0, \cdots, y_n),$$
}
where 
\begin{equation}\label{e-sum-def}
w_n(x_0, \cdots, x_n):=\sum_{i=0}^n p^ix_i^{p^{n-i}}=
x_0^{p^n}+px_1^{p^{n-1}}+\cdots+p^nx_n.
\end{equation}
\item 
For any ring $A$ and any elements 
$\alpha:=(a_0, a_1, \cdots), \beta:=(b_0, b_1, \cdots) \in W(A)$, 
it holds that 
$$\alpha+\beta=(S_0(a_0, b_0), S_1(a_0, b_0, a_1, b_1), \cdots).$$
\end{enumerate}
\end{lem}

\begin{proof} 
For a proof, see \cite[Ch. II, \S 6, disscussion after Theorem 6]{Ser79}. 
\end{proof}

\begin{lem}\label{l-sum-homog}
We use the same notation as in Lemma \ref{l-sum-def}. 
We equip $\Z[x_0, y_0, x_1, y_1, \cdots]$ 
with the graded $\Z$-algebra structure such that 
\begin{enumerate}
\item 
any element $n$ of $\Z$ is an homogeneous element of degree zero, and  
\item 
all $x_0, y_0, x_1, y_1, \cdots$ are homogeneous elements 
with $\deg x_n=\deg y_n=p^n$. 
\end{enumerate}
Then, for any non-negative integer $n$, 
$S_n(x_0, y_0, \cdots, x_n, y_n)$ is a homogeneous element of 
$\Z[x_0, y_0, x_1, y_1, \cdots]$ whose degree is equal to $p^n$. 
\end{lem}

\begin{proof}
For a proof, see \cite[(17.1.18)]{Haz12}. 
\end{proof}

\subsection{De Rham--Witt complex}

In this subsection, we summarise some of results established by \cite{Ill79}.

Let $k$ be a perfect field of characteristic $p>0$ and 
let $X$ be a scheme of finite type over $k$. 
We have de Rham-Witt sheaves $W_n\Omega_X^r$ on $X$, 
which are $W_n\MO_X$-modules. 
De Rham-Witt sheaves are equipped with 
$W\MO_X$-module homomorphisms 
(cf. \cite[Ch. I, D\'efinition 1.1(V2'), Th\'eor\`eme 1.3, Th\'eor\`eme 2.17]{Ill79}): 
\[
F:W_{n+1}\Omega_X^r \to (F_X)_*(W_n\Omega_X^r),
\]
\[
V:(F_X)_*(W_n\Omega_X^r) \to W_{n+1}\Omega_X^r,
\]
\[
R:W_{n+1}\Omega_X^r \to W_n\Omega_X^r.
\]
Note that the above map $F$ is $W_n\MO_X$-linear because 
$F: \bigoplus_{r \geq 0} W_{n+1}\Omega_X^r \to  \bigoplus_{r \geq 0}W_n\Omega_X^r$ is a homorphism of graded algebras \cite[Ch. I, Th\'eor\`eme 2.17]{Ill79}. 
Also $V$ is $W_n\MO_X$-linear by \cite[Ch. I, (1.15.2)]{Ill79}. 
Furthermore, $R$ is surjective \cite[Ch. I, Th\'eor\`eme 1.3]{Ill79}. 
Then we set $W\Omega_X^r:=\varprojlim_n W_n\Omega_X^r$, 
where we consider $\{W_n\Omega_X^r\}_{n \in \Z_{>0}}$ as a projective system via $R$. 
If $X$ is an $N$-dimensional smooth variety over $k$, 
then it is known that $W_nX$ is Cohen--Macaulay and $W_n\Omega_X^N$ is 
a dualising sheaf of $W_nX$ \cite[Theorem 4.1]{Eke84}.

\begin{lem}\label{l-DRW-coherent}
Let $k$ be a perfect field of characteristic $p>0$ and 
let $X$ be a scheme of finite type over $k$. 
Then $W_n\Omega_X^r$ is a coherent $W_n\MO_X$-module 
for any $n \in \Z_{>0}$ and any $r \in \Z_{\geq 0}$. 
In particular, also $(F_X^e)_*W_n\Omega_X^r$ is a coherent $W_n\MO_X$-module 
for any $e \in \Z_{\geq 0}$. 
\end{lem}

\begin{proof}
By \cite[Ch. I, Proposition 1.13.1]{Ill79}, 
each $W_n\Omega_X^r$ is a quasi-coherent $W_n\MO_X$-module. 
Assuming $X=\Spec\,A$, it suffices to show that 
$W_n\Omega_A^r=\Gamma(X, W_n\Omega_X^r)$ is a finitely generated $W_n(A)$-module. 
By \cite[Th\'eor\`eme 1.3]{Ill79}, 
there is a surjecitve $W_n(A)$-module homomorphism 
$$\Omega_{W_n(A)/\Z}^r \to W_n\Omega^r_A.$$
Hence, it is enough to prove that $\Omega_{W_n(A)/\Z}^1$ is 
a finitely generated $W_n(A)$-module. 
By \cite[Theorem 25.1]{Mat89}, 
we have an exact sequence:  
$$\Omega^1_{W_n(k)/\Z}\otimes_{W_n(k)} W_n(A) \to \Omega^1_{W_n(A)/\Z} \to \Omega^1_{W_n(A)/W_n(k)} \to 0.$$
Since $\Omega^1_{W_n(k)/\Z}=0$ (cf. \cite[Ch. I, Lemma 1.7]{Ill79}), 
we get an isomorphism: 
$$\Omega^1_{W_n(A)/\Z} \xrightarrow{\simeq} \Omega^1_{W_n(A)/W_n(k)}.$$
Since $W_n(A)$ is a finitely generated $W_n(k)$-algebra (\cite[Corollaire 1.5.7]{Ill79}), 
$\Omega^1_{W_n(A)/W_n(k)}$ is a finitely generated $W_n(A)$-module, 
hence so is $\Omega^1_{W_n(A)/\Z}$. 
\end{proof}

\begin{dfn}\label{d-gr-def}
Let $k$ be a perfect field of characteristic $p>0$ and 
let $X$ be a smooth variety over $k$. 
Fix $n \in \Z_{>0}$ and $r \in \Z_{\geq 0}$. 
We define a $W_{n+1}\MO_X$-module ${\rm gr}^n W\Omega_X^r$  by 
\[
{\rm gr}^n W\Omega_X^r := \Ker(R: W_{n+1}\Omega_X^r \to W_n\Omega_X^r). 
\]
It follows from Lemma \ref{l-DRW-coherent} that ${\rm gr}^n W\Omega_X^r$ is a coherent $W_{n+1}\MO_X$-module. 
\end{dfn}

\begin{rem}\label{r-gr-def}
Note that the $W_{n+1}\MO_X$-module $(F_X)_*({\rm gr}^n W\Omega_X^r)$ can be naturally considered as an $\MO_X$-module. 
Indeed, for any $a \in W_n\MO_X$ and $\zeta \in {\rm gr}^n W\Omega_X^r$, it holds that 
\[
(Va) \cdot (F_X)_*(\zeta) = (F_X)_*((FVa) \cdot \zeta)=(F_X)_*(pa \zeta)=0,
\]
where $(F_X)_*(\xi) \in (F_X)_*({\rm gr}^n W\Omega_X^r)$ denotes the element corresponding to 
$\xi \in {\rm gr}^n W\Omega_X^r$.  
Note that the last equality follows from the equation 
\[
{\rm gr}^n W\Omega_X^r  = \Ker(p: W_{n+1}\Omega_X^r \to W_{n+1}\Omega_X^r) 
\]
which is  guaranteed by \cite[Ch. I, Proposition 3.4]{Ill79}.
\end{rem}

\begin{prop}\label{p-DRW-proj-system}
Let $k$ be a perfect field of characteristic $p>0$ and 
let $X$ be a smooth variety over $k$. 
Fix $n \in \Z_{>0}$ and $r \in \Z_{\geq 0}$. 
Then the following hold. 
\begin{enumerate}
\item 
$(F_X)_*({\rm gr}^n W\Omega_X^r)$ is a coherent locally free $\MO_X$-module. 
\item 
There exists an exact sequences of $\MO_X$-modules: 
\[
0 \to \frac{(F_X^{n+1})_*\Omega_X^r}{(F_X)_*B_n\Omega_X^r} 
\xrightarrow{V^n} 
(F_X)_*({\rm gr}^n W\Omega_X^r) \xrightarrow{\beta} 
\frac{(F_X^{n+1})_*\Omega_X^{r-1}}{(F_X)_*Z_n\Omega_X^{r-1}} \to 0, 
\]
where the first $\MO_X$-module homomorphism  
\[
V^n :  (F_X^{n+1})_*\Omega_X^r / (F_X)_*B_n\Omega_X^r
\to 
(F_X)_*({\rm gr}^n W\Omega_X^r)
\]
is induced by $V^n : (F_X^n)_*\Omega_X^r \to W_{n+1}\Omega_X^r$ and 
the second $\MO_X$-module homomorphism 
\[
\beta: (F_X)_*({\rm gr}^n W\Omega_X^r) \to (F_X^{n+1})_*\Omega_X^{r-1} / (F_X)_*Z_n\Omega_X^{r-1}
\] 
satisfies 
\[
\beta ( V^n x + dV^n y ) = y + (F_X)_*Z_n\Omega_X^{r-1}
\] 
for $x \in (F_X^{n+1})_*\Omega_X^r$ and $y \in (F_X^{n+1})_*\Omega_X^{r-1}$. 
\end{enumerate}
\end{prop}

\begin{proof}
By Definition \ref{d-gr-def} and Remark \ref{r-gr-def}, 
$(F_X)_*{\rm gr}^n W\Omega_X^r$ is a coherent $\MO_X$-module. 
Since an extension of coherent locally free $\MO_X$-modules is a 
coherent locally free $\MO_X$-module, it follows from Proposition \ref{p-BZ}(3) that (2) implies (1). Hence it suffices to prove (2). 
The assertion (2) holds by 
the horizontal exact sequence in \cite[Ch. I, Corollaire 3.9]{Ill79}. 
Note that the $\MO_X$-module ${\rm gr}^n W\Omega_X^r$ in \cite[Ch. I, Corollaire 3.9]{Ill79} 
is nothing but our $(F_X)_*({\rm gr}^n W\Omega_X^r)$. 
Remark that there are the following typographic errors in \cite[Ch. I, Corollaire 3.9]{Ill79}. 
All of $B_n\Omega_X^i, B_{n+1}\Omega_X^i, Z_n\Omega_X^{i-1}$, and $Z_{n+1}\Omega_X^{i-1}$ should be 
replaced by 
$F_*B_n\Omega_X^i, F_*B_{n+1}\Omega_X^i, F_*Z_n\Omega_X^{i-1}$, and $F_*Z_{n+1}\Omega_X^{i-1}$, 
respectively. 
\end{proof}

\begin{thm}\label{t-descent}
Let $k$ be a perfect field of characteristic $p>0$ and 
let $X$ be an $N$-dimensional smooth variety over $k$. 
Then, for any positive integer $e$, the Frobenius homomorphism 
$$F^e:W\Omega_X^N \to (F_X^e)_*(W\Omega_X^N)$$
is a $W\MO_X$-module isomorphism. 
\end{thm}

\begin{proof}
The problem is reduced to the case when $e=1$. 
We may assume that $X$ is affine. 
By abuse of notation, we denote $\Gamma(X, M)$ simply by $M$ for a sheaf $M$ on $X$. 
Recall that $F$ can be written as follows: 
\begin{eqnarray*}
F: W\Omega_X^N = \varprojlim_n W_{n+1}\Omega_X^N 
&\to&  \varprojlim_n (F_X)_*W_{n}\Omega_X^N 
= (F_X)_*W\Omega_X^N\\
(a_{n+1})_n &\mapsto& (F(a_{n+1}))_n, 
\end{eqnarray*}
where $a_{n+1} \in  W_{n+1}\Omega_X^N$ and $F(a_{n+1}) \in (F_X)_*W_n\Omega_X^N$.

It follows from \cite[Ch. I, Proposition 3.7(b)]{Ill79} that 
$F: W_{n+1}\Omega_X^N \to (F_X)_*W_n\Omega_X^N$ is surjective for any $n \in \Z_{>0}$. 
For the kernel $K_n$ of $F: W_{n+1}\Omega_X^N \to (F_X)_*W_n\Omega_X^N$, 
we have an exact sequence of projective systems indexed by $\Z_{>0}$: 
\[
0 \to \{ K_n\}_n \to \left\{ W_{n+1}\Omega_X^N \right\}_n \xrightarrow{F} 
\left\{ (F_X)_*W_n\Omega_X^N \right\}_n \to 0. 
\]
It suffices to show that $\varprojlim_n K_n =0$ and 
the projective system $\{K_n\}_n$ satisfies the Mittag--Leffler condition.

It is enough to prove that $R : K_{n+1} \to K_n$ is zero. 
Pick $a_{n+1} \in K_{n+1}$. 
We have $a_{n+1} \in W_{n+1}\Omega_X^N$ and $F(a_{n+1})=0$. 
By $F(a_{n+1})=0$, we obtain $p a_{n+1} =VF(a_{n+1})=0$ in $W_{n+1}\Omega_X^N$. 
It follows from [Ill79, Ch. I, Proposition 3.4] that $a_n = R(a_{n+1})=0$. 
\end{proof}

\begin{lem}\label{l-DRW-inje}
Let $k$ be a perfect field of characteristic $p>0$ and 
let $X$ be a  smooth variety over $k$. 
Then, for any $n \in \Z_{>0}$, $r \in \Z_{\geq 0}$, and non-empty open subset $U$ of $X$, 
the induced map 
$$\Gamma(U, W_n\Omega_X^r) \to (W_n\Omega_X^r)_{\xi}$$
is injective, where $(W_n\Omega_X^r)_{\xi}$ denotes the stalk of $W_n\Omega_X^r$ 
at the generic point $\xi$ of $X$. 
\end{lem}

\begin{proof}
We may assume that $U$ is an affine open subset of $X$. 
By Definition \ref{d-gr-def}, we have an exact sequence 
$$0 \to {\rm gr}^n W\Omega_X^r \to W_{n+1}\Omega_X^r \to W_n\Omega_X^r \to 0,$$
where $(F_X)_*{\rm gr}^n W\Omega_X^r$ is a coherent locally free $\MO_X$-module 
(Proposition \ref{p-DRW-proj-system}).  
Thus, the assertion follows from the snake lemma and induction on $n$. 
\end{proof}


\section{Witt divisorial sheaves}\label{s-Witt-div-sh}

\subsection{Basic properties}

In this subsection, we introduce Witt divisorial sheaves: 
$W\MO_X(D)$ and $W_n\MO_X(D)$ (Definition \ref{d-div-def}). 
We also establish some fundamental properties such as 
coherence of $W_n\MO_X(D)$ (Proposition \ref{p-div-coherent}) and 
invariance under linear equivalence (Lemma \ref{l-linear-eq}).

\begin{dfn}\label{d-div-def}
Let $X$ be an integral normal noetherian $\F_p$-scheme. 
Let $D$ be an $\R$-divisor on $X$. 
Then we define the subpresheaf $W\MO_X(D)$ of the constant sheaf $W(K(X))$ on $X$ 
by 
{\small 
$$\Gamma(U, W\MO_X(D)):=\{(\varphi_0, \varphi_1, \cdots) \in W(K(X))\,|\,
\left({\rm div}(\varphi_n)+p^nD\right)|_U \geq 0\}$$
}
for any open subset $U$ of $X$, 
where ${\rm div}(\varphi_n)$ denotes the principal divisor associated to $\varphi_n$. 
By definition, $W\MO_X(D)$ is a subsheaf of $W(K(X))$ (cf. Remark \ref{r-div-def}). 
We call $W\MO_X(D)$ the {\em Witt divisorial sheaf associated to} $D$. 
We define the subsheaf $W_n\MO_X(D)$ of the constant sheaf $W_n(K(X))$ 
in the same way.  
In particular, we have 
\[
\Gamma(U, W_1\MO_X(D)) = 
\Gamma(U, \MO_X(D)) = 
\{ \varphi \in K(X) \,|\,
\left({\rm div}(\varphi)+D\right)|_U \geq 0\}. 
\]
\end{dfn}

\begin{rem}\label{r-div-def}
We use notation as in Definition \ref{d-div-def}. If we identify $W(K(X))$ with the infinite direct product $\prod_{n=0}^{\infty}K(X)$ as sets, then it follows by definition that as subsets of $W(K(X))=\prod_{n=0}^{\infty}K(X)$, we obtain an equation $$\Gamma(U, W\MO_X(D))=\prod_{n=0}^{\infty} \Gamma(U, \MO_X(p^nD)).$$
\end{rem}

\begin{rem}\label{r-perturb}
We use notation as in Definition \ref{d-div-def}. 
Fix a positive integer $n$. 
By Remark \ref{r-div-def}, we obtain an equation of subsets of 
$W_n(K(X))=\prod_{m=0}^{n-1}K(X)$: 
$$\Gamma(U, W_n\MO_X(D))=\prod_{m=0}^{n-1} \Gamma(U, \MO_X(p^mD))
=\prod_{m=0}^{n-1} \Gamma(U, \MO_X(\llcorner p^mD\lrcorner)).$$
Thus, if $E$ is an effective $\R$-divisor on $X$, 
then there exists a positive real number $\epsilon$ such that 
the equation 
$$W_n\MO_X(D)=W_n\MO_X(D')$$
holds as subsheaves of the constant sheaf $W_n(K(X))$ for any $\R$-divisor $D'$ 
satisfying $D \leq D' \leq D+\epsilon E$. 
\end{rem}

\begin{rem}\label{r-BK}
Let $K$ be a discrete valuation field. 
Brylinski--Kato introduced the following subgroup of $W_n(K)$ \cite[Definition 3.1]{Kat89} 
(strictly speaking, \cite{Kat89} assumes $K$ to be henselian): 
\[
{\rm fil}_r(W_n(K)) = \{(b_0, ..., b_{n-1}) \in W_n(K)\,|\, 
{\rm ord}_K(b_i) \geq -\frac{rp^i}{p^{n-1}} \text{ for all }i \}.
\]
This subgroup can be considered as a stalk of a Witt divisorial sheaf as follows.

Let $X$ be an integral normal noetherian $\F_p$-scheme and let $D$ be a $\Z$-divisor on $X$. 
Set $K:=K(X)$. 
Fix an irreducible component $D'$ of $D$. 
Let $\xi$ be the generic point of $D'$ and 
let $r$ be the coefficient of $D'$ in $D$, i.e. $D=rD'+\cdots$. 
We consider $K$ as a discrete valuation field with respect to the valuation corresponding to $\MO_{X, \xi}$. 
Then the equation 
\[
W_n\MO_X\left( 
(1/p^{n-1})D\right)_{\xi} = {\rm fil}_r (W_n(K)) 
\]
holds. 
\end{rem}

\begin{lem}\label{l-div-sub}
Let $X$ be an integral normal noetherian $\F_p$-scheme. 
Let $D$ be an $\R$-divisor on $X$. 
Then the following hold. 
\begin{enumerate}
\item 
$W\MO_X(D)$ is a sheaf of $W\MO_X$-submodules of $W(K(X))$, 
i.e. for any open subset $U$ of $X$, 
$\Gamma(U, W\MO_X(D))$ is a $W\MO_X(U)$-submodule of $W(K(X))$. 
\item 
For any positive integer $n$, 
$W_n\MO_X(D)$ is a sheaf of $W_n\MO_X$-submodules of $W_n(K(X))$, 
i.e. for any open subset $U$ of $X$, 
$\Gamma(U, W_n\MO_X(D))$ is a $W_n\MO_X(U)$-submodule of $W_n(K(X))$. 
Furthermore, $W_n\MO_X(D)$ is a quasi-coherent $W_n\MO_X$-module. 
\end{enumerate}
\end{lem}

\begin{proof}
Let us show (1). 
We may assume that $X=U$ and $X$ is affine, say $X=U=\Spec\,A$. 

Fix $\varphi, \psi \in \Gamma(X, W\MO_X(D))$. 
Let us prove $\varphi+\psi \in \Gamma(X, W\MO_X(D))$. 
We can write 
$$\varphi=(\varphi_0, \varphi_1, \cdots), \quad \psi=(\psi_0, \psi_1, \cdots)$$
for some $\varphi_n, \psi_n \in \Gamma(X, \MO_X(p^nD))$. 
By Lemma \ref{l-sum-def}, 
we have that 
$$\varphi+\psi=(S_0(\varphi_0, \psi_0), S_1(\varphi_0, \psi_0, \varphi_1, \psi_1), 
\cdots)$$
for some polynomials 
$$S_n(x_0, y_0, \cdots, x_n, y_n) \in \Z[x_0, y_0, \cdots, x_n, y_n]$$
satisfying the properties listed in Lemma \ref{l-sum-def}. 
We equip the polynomial ring $\Z[x_0, y_0, x_1, y_1, \cdots]$ with the structure of graded $\Z$-algebra defined in the statement of Lemma \ref{l-sum-homog}, 
i.e. we consider $\Z[x_0, y_0, x_1, y_1, \cdots]$ as a weighted  polynomial ring 
such that $\deg x_i=\deg y_i=p^i$. 
Pick a monomial 
$x_0^{a_0}y_0^{b_0} \cdots x_n^{a_n}y_n^{b_n}$ appearing in 
the monomial decomposition of $S_n(x_0, y_0, \cdots, x_n, y_n)$. 
Since $S_n$ is homogeneous of degree $p^n$ (Lemma \ref{l-sum-homog}), 
it holds that 
$$\sum_{i=0}^n p^i(a_i+b_i)=p^n.$$
By 
$${\rm div} (\varphi_i) \geq -p^iD\qquad {\rm and}\quad \quad {\rm div} (\psi_i) \geq -p^iD,$$
we have that 
\begin{eqnarray*}
{\rm div} (\varphi_0^{a_0}\psi_0^{b_0} \cdots \varphi_n^{a_n}\psi_n^{b_n})
&=&\sum_{i=0}^n (a_i{\rm div} (\varphi_i)+b_i{\rm div} (\psi_i))\\
&\geq &-\sum_{i=0}^n p^i(a_i+b_i)D\\
&=& -p^nD.
\end{eqnarray*}
In other words, we obtain 
$\varphi_0^{a_0}\psi_0^{b_0} \cdots \varphi_n^{a_n}\psi_n^{b_n} 
\in \Gamma(X, \MO_X(p^nD))$. 
Therefore, it holds that $\varphi+\psi \in \Gamma(X, W\MO_X(D))$.

Fix $a=(a_0, a_1, \cdots) \in W(A)$ and $\varphi=(\varphi_0, \varphi_1, \cdots) \in \Gamma(X, W\MO_X(D))$. 
Let us show $a\varphi \in \Gamma(X, W\MO_X(D))$. 
We can write 
$$a\varphi=\left(\sum_{m=0}^{\infty} V^m({\underline{a_m}})\right) \cdot (\varphi_0, \varphi_1, \cdots).$$
Thus, for any $b \in A$ and any $m \in \Z_{\geq 0}$, it suffices to prove that 
$$(V^m(\underline{b})) \cdot (\varphi_0, \varphi_1, \cdots) \in \Gamma(X, W\MO_X(D)).$$
We have that 
\begin{eqnarray*}
(V^m(\underline{b})) \cdot (\varphi_0, \varphi_1, \cdots)
&=&V^m(\underline{b} \cdot (\varphi_0^{p^m}, \varphi_1^{p^m}, \cdots))\\
&=&V^m(b\varphi_0^{p^m}, b^p\varphi_1^{p^m}, \cdots)\\
&=:&(\psi_0, \psi_1, \cdots). 
\end{eqnarray*}
In particular, we get  $\psi_0=\cdots=\psi_{m-1}=0$. 
Since we have ${\rm div}(\varphi_{\ell}) \geq -p^{\ell}D$ for any $\ell$, 
the following equation holds for any $n \geq m$: 
$$
{\rm div}(\psi_n) = {\rm div}(b^{p^{n-m}}\varphi_{n-m}^{p^m}) 
\geq {\rm div}(\varphi_{n-m}^{p^m})$$
$$=p^m {\rm div}(\varphi_{n-m}) \geq p^m\cdot (-p^{n-m}D)=-p^nD.$$
Therefore, we get $a\varphi \in \Gamma(X, W\MO_X(D))$. 
Thus, (1) holds. 

Let us show (2). 
By the same argument as in (1), we have that $W_n\MO_X(D)$ is 
a sheaf of $W_n\MO_X$-submodules of $W_n(K(X))$. 
What is remaining is to prove that $W_n\MO_X(D)$ is 
a quasi-coherent $W_n\MO_X$-module. 
We may assume that $X=\Spec\,A$. 
Take $f \in A \setminus \{0\}$ 
 and let $\underline{f} \in W_n(A)$ be the Teichm\"{u}ller lift  of $f$. 
It suffices to show that the natural $W_n(A)_{\underline{f}}$-module homomorphism 
$$\theta:\Gamma(X, W_n\MO_X(D))_{\underline{f}} \to \Gamma(D(f), W_n\MO_X(D))$$
is bijective. 
Since $\theta$ is compatible with the inclusions to $W_n(K(X))$, 
it holds that $\theta$ is injective. 
Let us prove that $\theta$ is surjective. 
Take 
$$\varphi=(\varphi_0, \cdots, \varphi_{n-1}) \in \Gamma(D(f), W_n\MO_X(D)),$$ 
where $({\rm div}(\varphi_i)+p^iD)|_{D(f)} \geq 0$. 
We can find $N \in \Z_{>0}$ such that 
$${\rm div}(f^N\varphi_i)+p^iD \geq 0$$
for any $i \in \{0, \cdots, n-1\}$. 
In particular, 
$$\underline{f}^N\varphi=(f^N\varphi_0, \cdots, f^{Np^{n-1}}\varphi_{n-1}) \in \Gamma(X, W_n\MO_X(D)).$$
Thus, $\theta$ is surjective. 
Hence, (2) holds. 
\end{proof}

\begin{rem}\label{r-div-codim2}
We use notation as in Definition \ref{d-div-def}. 
Take an open subset $X'$ of $X$ such that 
$X'$ contains all the points of codimension one in $X$. 
Let $j:X' \to X$ be the induced open immersion. 
Then it holds that both 
$W\MO_X(D) \to j_*(W\MO_X(D)|_{X'})$ and $W_n\MO_X(D) \to j_*(W_n\MO_X(D)|_{X'})$ 
are isomorphisms for any $n \in \Z_{>0}$ (cf. Remark \ref{r-div-def}). 
\end{rem}

\begin{prop}\label{p-div-induction}
Let $X$ be an integral normal noetherian $\F_p$-scheme and 
let $D$ be an $\R$-divisor on $X$. 
Then there are exact sequences of $W\MO_X$-module homomorphisms: 
\begin{equation}\label{e-div-induction1}
0 \to (F_X^n)_*(W\MO_X(p^nD)) \xrightarrow{V^n} W\MO_X(D) \to W_n\MO_X(D) \to 0
\end{equation}
{\small 
\begin{equation}\label{e-div-induction2}
0 \to (F_X^n)_*(W_m\MO_X(p^nD)) \xrightarrow{V^n} W_{n+m}\MO_X(D) \to W_n\MO_X(D) \to 0
\end{equation}
}
for any positive integers $n$ and $m$. 
\end{prop}

\begin{proof}
Let us prove that 
the first sequence (\ref{e-div-induction1}) is exact. 
This is a subcomplex of the exact sequence (cf. \cite[Ch. 0, (1.1.6)]{Ill79}): 
$$0 \to (F_X^n)_*W(K(X)) \xrightarrow{V^n} W(K(X)) \to W_n(K(X)) \to 0.$$
In particular, $V^n:(F_X^n)_*(W\MO_X(p^nD)) \to W\MO_X(D)$ is injective. 
By construction, 
the latter homomorphism $W\MO_X(D) \to W_n\MO_X(D)$ is surjective. 
Let us prove the exactness on the middle term $W\MO_X(D)$. 
Take an element
$$\varphi =(\varphi_0, \varphi_1, \cdots) \in \Gamma(X, W\MO_X(D))$$
whose image to $W_n\MO_X(D)$ is zero, 
i.e. $\varphi_0=\cdots=\varphi_{n-1}=0$. 
We can check directly from Definition \ref{d-div-def} that 
$$(\varphi_n, \varphi_{n+1}, \cdots) \in \Gamma(X, W_n\MO_X(p^nD)),$$
hence the sequence (\ref{e-div-induction1}) is exact. 
It holds by the same argument that also (\ref{e-div-induction2}) is exact. 
\end{proof}

\begin{prop}\label{p-div-coherent}
Let $X$ be an integral normal $F$-finite noetherian $\F_p$-scheme. 
Let $D$ be an $\R$-divisor on $X$. 
Then, for any positive integer $n$, 
$W_n\MO_X(D)$ is a coherent $W_n\MO_X$-module. 
\end{prop}

\begin{proof}
We may assume that $X$ is affine. 
We prove the assertion by induction on $n$. 
Let us first treat the case when $n=1$, although this case is well known to experts. 
By Lemma \ref{l-div-sub}(2), 
$\MO_X(D)$ is a quasi-coherent $\MO_X$-module. 
Fix an effective Cartier divisor $E$ such that $E-D$ is effective and ${\rm div}(f) =E$ for some $f \in \Gamma(X, \MO_X)$. 
We have an isomorphism of $\MO_X$-modules: 
\[
\MO_X(D) \xrightarrow{\simeq, \times f} \MO_X(-E+D) = \MO_X(-(E-D)). 
\]
Then $\MO_X(-(E-D))$ is a coherent $\MO_X$-module, because this is 
a quasi-coherent $\MO_X$-submodule of $\MO_X$. 
This completes the proof for the case when $n=1$. 

By Lemma \ref{l-div-sub}(2), 
$W_n\MO_X(D)$ is a quasi-coherent $W_n\MO_X$-module. 
Since $X$ is $F$-finite, $(F^n_X)_*(M)$ is a coherent $W_n\MO_X$-module 
for any coherent $W_n\MO_X$-module $M$. 
Thus the assertion follows from the induction hypothesis 
and the exact sequence (\ref{e-div-induction2}) of Lemma \ref{p-div-induction}. 
\end{proof}

\begin{lem}\label{l-linear-eq}
Let $X$ be an integral normal noetherian $\F_p$-scheme. 
Assume that $\R$-divisors $D_1$ and $D_2$ on $X$ are $\Z$-linearly equivalent i.e. 
$D_2=D_1+{\rm div}(\varphi)$ for some $\varphi \in K(X)^{\times}$. 
Then 
$W\MO_X(D_1)$ and $W\MO_X(D_2)$ are isomorphic as $W\MO_X$-modules, 
and 
$W_n\MO_X(D_1)$ and $W_n\MO_X(D_2)$ are isomorphic as $W_n\MO_X$-modules 
for any $n \in \Z_{>0}$. 
\end{lem}

\begin{proof}
We have an $\MO_X$-module isomorphism: 
$$\times \varphi^{p^n}: \MO_X(p^nD_1) \xrightarrow{\simeq} \MO_X(p^nD_2), \quad f \mapsto \varphi^{p^n} f.$$
Therefore, the $W\MO_X$-module homomorphism 
$$\times \underline{\varphi}: W(K(X)) \to W(K(X))$$
induces a $W\MO_X$-module isomorphism $W\MO_X(D_1) \simeq W\MO_X(D_2)$. 
The same argument implies $W_n\MO_X(D_1) \simeq W_n\MO_X(D_2)$ 
for any $n \in \Z_{>0}$. 
\end{proof}

\subsection{Teichm\"{u}ller lifts of invertible sheaves}

In this subsection, we introduce Teichm\"uller lifts of invertible sheaves 
(Definition \ref{d-teich-lift}). 
This notion will turn out to be a special case of Witt divisorial sheaves 
for normal varieties (Proposition \ref{p-div-Teich}).

\begin{dfn}\label{d-teich-lift}
Let $X$ be an $\F_p$-scheme and let $L$ be an invertible $\MO_X$-module. 
Then $L$ can be represented by 
 $\{(U_i)_{i \in I}, (f_{ij})_{(i, j) \in I \times I}\}$ satisfying 
the following properties:  
\begin{itemize}
\item 
$U_i$ is an open affine subset of $X$ for any $i \in I$. 
\item 
$X=\bigcup_{i \in I} U_i$. 
\item 
$f_{ij} \in \Gamma(U_{ij}, \MO_X)^{\times}$ for any $i, j \in I$,  
where $U_{ij} := U_i \cap U_j$. 
\item 
The equation $(f_{ij}|_{U_{ijk}})(f_{jk}|_{U_{ijk}})(f_{ki}|_{U_{ijk}})=1$ 
holds in $\Gamma(U_{ijk}, \MO_X)^{\times}$,  
where $U_{ijk} := U_i \cap U_j \cap U_k$. 
\end{itemize}
We define the {\em Teichm\"uller lift of} $L$, denoted by $\underline{L}$, 
as the invertible $W\MO_X$-module defined by 
$\{(U_i)_{i \in I}, (\underline{f_{ij}})_{(i, j) \in I \times I}\}$.  
Note that $\underline L$ coincides with the invertible $W\MO_X$-module 
corresponding to the image of $L$ by the group homomorphism  
$$\check{H}^1(X, \MO_X^{\times}) \to \check{H}^1(X, W\MO_X^{\times})$$
that is induced by the Teichm\"{u}ller lift 
$\MO_X^{\times} \to W\MO_X^{\times},\,\, f \mapsto \underline f$. 
In particular, the isomorphism class of $\underline L$ does not depend on the choice of 
representation $\{(U_i)_{i \in I}, (f_{ij})_{(i, j) \in I \times I}\}$. 
For any positive integer $n$, 
we define the invertible $W_n\MO_X$-module $\underline{L}_{\leq n}$ in the same way. 
\end{dfn}

\begin{rem}\label{r-teich-lift}
By definition, we obtain 
$\underline{L} \otimes_{W\MO_X} W_n\MO_X \simeq \underline{L}_{\leq n}$ 
and $\underline{L}_{\leq 1} =L$. 
\end{rem}

\begin{prop}\label{p-div-Teich}
Let $X$ be an integral normal noetherian $\F_p$-scheme. 
Let $D$ be a Cartier divisor on $X$. 
Then the following hold. 
\begin{enumerate}
\item 
The Witt divisorial sheaf $W\MO_X(D)$ 
is isomorphic to the Teichm\"uller lift $\underline{\MO_X(D)}$ of $\MO_X(D)$ 
as $W\MO_X$-modules. 
\item 
The Witt divisorial sheaf $W_n\MO_X(D)$ 
is isomorphic to the Teichm\"uller lift $\underline{\MO_X(D)}_{\leq n}$ of $\MO_X(D)$ 
as $W_n\MO_X$-modules. 
\end{enumerate}
\end{prop}

\begin{proof}
Let us show (1). 
By definition of Cartier divisors, $D$ is represented by $\{(U_i, f_i)\}$, 
where $X=\bigcup_{i \in I} U_i$ is an affine open cover of $X$, 
$f_i \in K(X)^{\times}$ and $f_i/f_j \in \Gamma(U_i \cap U_j, \MO_X^{\times})$. 
Then, for any open subset $V$ of $X$, we have that 
$$\Gamma(V, W\MO_X(D))=
\bigcap_{i \in I}\{(\varphi_0, \varphi_1, \cdots) \in W(K(X))\,|\, 
{\rm div}( f_i^{p^n}\varphi_n)|_{V \cap U_i} \geq 0\}.$$
Therefore, we get an equation 
$$W\MO_X(D)|_{U_i}=\underline{f_i^{-1}}W\MO_{U_i}$$
as subsheaves of the constant sheaf $W(K(X))$, 
where the right hand side is isomorphic to $W\MO_{U_i}$. 
In particular, $W\MO_X(D)$ is the invertible $W\MO_X$-module represented by 
$\{(U_i)_{i \in I}, (\underline{f_i^{-1}f_j})_{(i, j) \in I \times I}\}$, 
which is isomorphic to the Teichm\"uller lift $\underline{\MO_X(D)}$. 
Thus (1) holds. 
We omit the proof of (2), as we can apply the same argument as in (1). 
\end{proof}

\subsection{Pullback of divisors}\label{ss-def-pullback-div}


In this subsection, we recall definition of pullbacks of $\R$-divisors 
and introduce a natural homomorphism 
\[
W\MO_Y(D) \to f_*W\MO_X(f^*D)
\]
for some cases.

\begin{nothing}[Definition of pullback]\label{d-pullback}
Let $f:X \to Y$ be a dominant morphism of integral normal excellent $\F_p$-schemes. 
Let $D$ be an $\R$-divisor on $Y$. 
Since $f$ is dominant, we have a field extension $\varphi:K(Y) \hookrightarrow K(X)$, 
which induces an injective ring homomorphism  
\begin{equation*}
W(\varphi):W(K(Y)) \hookrightarrow W(K(X)). 
\end{equation*}
In the following, we shall define $f^*D$ and a $W\MO_Y$-module homomorphism 
$$W\MO_Y(D) \to f_*W\MO_X(f^*D)$$
compatible with $W(\varphi)$, under the assumption that 
\begin{enumerate}
\item $D$ is $\R$-Cartier, or 
\item $f$ is separated and quasi-finite. 
\end{enumerate}

(1) First we assume that $D$ is $\R$-Cartier. 
Then $f^*D$ is defined as in Subsection \ref{ss-notation}(\ref{ss-notation-9}), hence we have 
$$\varphi':\MO_Y(p^mD) \to \MO_X(p^mf^*D)$$
that is compatible with $\varphi$. Thus, $W(\varphi)$ induces 
$W\MO_Y$-module homomorphism: 
$$\alpha: W\MO_Y(D) \to W\MO_X(f^*D).$$ 

(2)' As a special case of (2), we treat the case when $f$ is a finite morphism. 
Let $X_0$ and $Y_0$ be the regular loci of $X$ and $Y$ respectively. 
Note that $X_0$ and $Y_0$ are open subsets of $X$ and $Y$ respectively, 
because $X$ and $Y$ are assumed to be excellent schemes. 
We set 
$$Y_1:=Y_0 \cap (Y \setminus f(X \setminus X_0)), \quad X_1:=f^{-1}(Y_1)$$
and let $f_1:X_1 \to Y_1$ be the induced morphism. 
We define $f^*D$ as the closure of the $\R$-divisor $f_1^*D|_{Y_1}$ on $X_1$. 
By (1), we have an $W\MO_{Y_1}$-module homomorphism 
$$\alpha: W\MO_Y(D)|_{Y_1} \to (f_1)_*(W\MO_X(f^*D)|_{X_1})=(f_*W\MO_X(f^*D))|_{Y_1}$$
compatible with $W(\varphi)$. 
By Remark \ref{r-div-codim2}, we get a $W\MO_Y$-module homomorphism  
$$\beta:=j_*\alpha: W\MO_Y(D) \to W\MO_X(f^*D)$$
compatible with $W(\varphi)$, 
where $j:Y_1 \hookrightarrow Y$ is the induced open immersion.

(2) Finally we treat the case when $f$ is a separated quasi-finite morphism. 
By \cite[Theorem 1.10.13]{Fu11}, 
there is a factorisation 
$$f:X \xrightarrow{i} Z \xrightarrow{g} Y$$
where $Z$ is an integral normal excellent $\F_p$-scheme, 
$i$ is an open immersion and $g$ is a finite morphism. 
We define $f^*D:=(g^*D)|_X$, 
where $g^*D$ is defined in (2)'. 
By (2)', we have an $W\MO_Y$-module homomorphism 
$$\beta: W\MO_Y(D) \to g_*W\MO_Z(g^*D),$$
which induces $W\MO_Y$-module homomorphisms 
$$\gamma: W\MO_Y(D) \xrightarrow{\beta} g_*W\MO_Z(g^*D) 
\xrightarrow{\delta} g_*i_*W\MO_X(g^*D|_X)=f_*W\MO_X(f^*D),$$
where $\delta$ is the restriction map. 

By the same argument as above, 
if one of (1) and (2) holds, then 
we obtain a $W_n\MO_Y$-module homomorphism 
$W_n\MO_Y(D) \to f_*W_n\MO_X(f^*D)$. 
\end{nothing}

\subsection{CM-ness of $W_n\MO_X(D)$}

In this subsection, we prove that $W_n\MO_X(D)$ is Cohen-Macaulay 
under suitable conditions (Theorem \ref{t-div-CM}). 
As an application, we obtain some vanishing of Ext sheaves (Proposition \ref{p-div-CM}).

\begin{thm}\label{t-div-CM}
Let $X$ be an integral regular noetherian $\F_p$-scheme. 
Let $D$ be an $\R$-divisor on $X$. 
Then, for any positive integer $n$, 
$W_n\MO_X(D)$ is a maximal Cohen-Macaulay $W_n\MO_X$-module. 
\end{thm}

\begin{proof}
Since the problem is local on $X$, we may assume that $X$ is affine. 
Fix a closed point $x$ of $X$. 
In what follows, $H^i_x(M)$ denotes the local cohomology of a sheaf $M$ of abelian groups 
in the sense of \cite[Definition in page 2]{Har67}. 
If $M$ is a coherent $W_n\MO_X$-module on $W_nX$, 
then $H^i_x(M)$ coincides with $H^i_{\m_n}(M_x)$ in the sense of \cite[Definition 3.5.2]{BH93}, where $M_x$ is the $\MO_{W_nX, x}$-module and $\m_n$ denotes the maximal ideal of the local ring $\MO_{W_nX, x}$. 

By $X=\Supp\,(W_n\MO_X(D))$, it is enough to show that $W_n\MO_X(D)_x$ is 
a Cohen-Macaulay $W_n\MO_{X, x}$-module \cite[Theorem 17.3(iii)]{Mat89}. 
We may assume that $X=\Spec\,A$, where $A$ is a local ring. 
Set $d:=\dim A$. 

Let us prove $H^i_x(W_n\MO_X(D)) =0$ for $i \neq d$ by induction on $n$. 
We first treat the case when $n=1$. 
Since $A$ is a regular local ring, 
it holds that $W_1\MO_X(D)=\MO_X(D) \simeq \MO_X$. 
Therefore, we obtain 
$H^i_x(W_1\MO_X(D)) \simeq H^i_x(\MO_{X})=0$ for $i \neq d$, 
because $X$ is Cohen--Macaulay. 

Assume $n \geq 1$. 
By Proposition \ref{p-div-induction}, 
we have an exact sequence of the sheaves of abelian groups: 
\[
0 \to \MO_X(p^{n+1}D) \to W_{n+1}\MO_X(D) \to W_n\MO_X(D) \to 0. 
\]
By the induction hypothesis, it holds that 
$H^i_x(\MO_{X}(p^{n+1}D))=H^i_x(W_n\MO_X(D))=0$ for $i \neq d$. 
Then it follows from the above exact sequence that $H_x^i(W_{n+1}\MO_X(D))=0$ for $i \neq d$. 
\end{proof}

\begin{prop}\label{p-div-CM}
Let $k$ be a perfect field of characteristic $p>0$ and 
let $X$ be an $N$-dimensional smooth variety over $k$. 
Let $D$ be an $\R$-divisor on $X$. 
Then, for any non-negative integer $e$ and positive integers $n$ and $m$ satisfying $n \geq m$, the isomorphism 
$$R\mathcal Hom_{W_n\MO_X}((F_{W_nX}^e)_*j_*W_m\MO_X(D), W_n\Omega_X^N)$$ 
$$\simeq (F_{W_nX}^e)_*j_*\mathcal Hom_{W_m\MO_X}(W_m\MO_X(D), W_m\Omega_X^N)$$
holds in the derived category of $W_n\MO_X$-modules, 
where $j:W_mX \to W_nX$ is the induced closed immersion.  
\end{prop}

\begin{proof}
The assertion follows from isomorphisms: 
\begin{eqnarray*}
&&R\mathcal Hom_{W_n\MO_X}((F_{W_nX}^e)_*j_*W_m\MO_X(D), W_n\Omega_X^N)\\
& \simeq &
(F_{W_nX}^e)_*j_* R\mathcal Hom_{W_m\MO_X}(W_m\MO_X(D), (F_{W_nX}^e \circ j)^!W_n\Omega_X^N)\\
& \simeq &
(F_{W_nX}^e)_*j_* R\mathcal Hom_{W_m\MO_X}(W_m\MO_X(D), W_m\Omega_X^N)\\
& \simeq &
(F_{W_nX}^e)_*j_* \mathcal Hom_{W_m\MO_X}(W_m\MO_X(D), W_m\Omega_X^N),
\end{eqnarray*}
where the first isomorphism holds by Grothendieck duality, 
the second isomorphism follows from \cite[Theorem 4.1]{Eke84}, 
and we obtain the third isomorphism by Theorem \ref{t-div-CM} and 
\cite[Corollary 3.5.11(a)]{BH93}. 
\end{proof}

\section{Vanishing theorems}

In this section, we prove the main theorems of this paper. 
In Subsection \ref{ss-Kodaira}, we prove a Kodaira vanishing theorem 
for $W\Omega_X^N$ (Theorem \ref{t-Kod-van}). 
Subsection \ref{ss4-split} and Subsection \ref{ss-WO-WO} are devoted to giving  
auxiliary results for Subsection \ref{ss-Kaw-Vie}. 
In Subsection \ref{ss-Kaw-Vie}, we establish Kawamata--Viehweg vanishing theorems 
for $W\Omega_X^N$ (Theorem \ref{t-KVV}, Theorem \ref{t-KVV-minus}).

\subsection{Kodaira vanishing}\label{ss-Kodaira}

The purpose of this subsection is to prove Theorem \ref{t-Kod-van}. 
We start by recalling a criterion for the vanishing of $R^q\varprojlim_n(-)$.

\begin{lem}\label{l-CR12}
Let $X$ be an $\F_p$-scheme and let $E=(E_n)_{n \geq 1}$ be a projective system of $W\MO_X$-modules. 
Consider the following three conditions. 
\begin{enumerate}
\item[(a)] $H^i(U, E_n)=0$ for any $i>0$, $n>0$, and affine open subset $U$ of $X$. 
\item[(b)] The projective system $\{H^0(U, E_n)\}_{n \geq 1}$ satisfies the Mittag--Leffler condition 
for any affine open subset $U$ of $X$. 
\item[(c)] 
The projective system $\{H^j(X, E_n)\}_{n \geq 1}$ satisfies the Mittag--Leffler condition for any $j\geq 0$. 
\end{enumerate}
Then the following hold. 
\begin{enumerate}
\item If $({\rm a})$ holds, then $R^q\varprojlim_n E_n =0$ for any $q\geq 2$. 
\item If $({\rm a})$ and $({\rm b})$ hold, then $R^q\varprojlim_n E_n =0$ for any $q \geq 1$. 
\item If $({\rm a}), ({\rm b}),$ and $({\rm c})$ hold, then the isomorphism 
\[
\varprojlim_n H^j(X, E_n) \simeq H^j(X, \varprojlim_n E_n)
\] 
holds for any $j \geq 0$. 
\end{enumerate}
\end{lem}

\begin{proof}
The assertions (1) and (2) follow from \cite[Lemma 1.5.1]{CR12}. 
Let us show (3). Set $E:=  \varprojlim_n E_n$. 
Let $\alpha : X \to \Spec\,\F_p$ be the induced morphism. 
We have isomorphisms 
\[
R\varprojlim_n \circ R\alpha_* E_n \simeq 
R\alpha_* \circ R\varprojlim_n E_n \simeq R\alpha_* \left( \varprojlim_n E_n \right) = R\alpha_*  E. 
\]
where the second isomorphism holds by (2). 
We obtain the following spectral sequence: 
\[
E_2^{q, j} := R^q\varprojlim_n ( H^j(X, E_n) ) \Rightarrow H^{q+j}( X, E) =:E^{q+j}.
\]
By (c), we have $E_2^{q, j}=0$ for any $q >0$ and $j \geq 0$. 
Then it holds that 
\[
\varprojlim_n ( H^j(X, E_n) ) =E_2^{0, j} \simeq E^j = H^{j}( X, E), 
\]
as required. 
\end{proof}

\begin{thm}\label{t-Kod-van}
Let $k$ be a perfect field of characteristic $p>0$. 
Let $f:X \to Y$ be a projective $k$-morphism 
from an $N$-dimensional smooth $k$-variety $X$ to 
a scheme $Y$ of finite type over $k$. 
Let $A$ be an $f$-ample Cartier divisor on $X$. 
Then the following hold. 
\begin{enumerate}
\item 
The isomorphism
\[
M \otimes_{W\MO_X} W\MO_X(sA) \simeq 
M \otimes_{W_{n'}\MO_X} W_{n'}\MO_X(sA)
\]
holds for any coherent $W_n\MO_X$-module $M$, $s \in \Z$, $n \in \Z$, and $n' \in \Z$ such that $n' \geq n > 0$. 
\item 
There exists a positive integer $s_1$ such that 
\[
R^if_*(W_n\Omega_X^N \otimes_{W\MO_X} W\MO_X(sA))  =0
\]
for any $i \in \Z_{>0}$, $n \in \Z_{>0}$, and $s \in \Z_{\geq s_1}$. 
\item 
There exists a positive integer $s_2$ such that 
\[
R^q\varprojlim_n f_*(W_n\Omega_X^N \otimes_{W\MO_X} W\MO_X(spA)) 
=0
\]
for any $q \in \Z_{>0}$ and $s \in \Z_{\geq s_2}$. 
\item 
The equation 
\[
R^if_*(W\Omega_X^N \otimes_{W\MO_X} W\MO_X(A)) =0
\]
holds for any $i \in \Z_{>0}$. 
\end{enumerate}
\end{thm}

\begin{proof}
The assertion (1) follows from the fact that 
$W\MO_X(sA)$ and $W_{n'}\MO_X(sA)$ are invertible sheaves whose transition functions 
are the same after taking the restriction to $W_n\MO_X$. 

Let us show (2). 
Applying $\mathcal Hom_{W_{n+1}\MO_X}(-, W_{n+1}\Omega_X^N)$ and its derived functor to 
\[
0 \to (F_X)_*W_n\MO_X \xrightarrow{V} W_{n+1}\MO_X \xrightarrow{R} \MO_X \to 0,
\]
it follows from Proposition \ref{p-div-CM} that we obtain the exact sequence: 
\[
0 \to \Omega_X^N \xrightarrow{R^*} W_{n+1}\Omega_X^N \xrightarrow{V^*} (F_X)_*(W_n\Omega_X^N) 
\to 0. 
\]
Then the assertion (2) follows from the Serre vanishing theorem and induction on $n$. 

Let us show (3). 
Set 
\[
E_n(s) := f_*(W_n\Omega_X^N \otimes_{W\MO_X} W\MO_X(spA)). 
\]
Assume that $Y$ is affine. 
By Lemma \ref{l-CR12}(2), it suffices to find $s_2 \in \Z_{>0}$ 
which satisfies the following conditions (a) and (b). 
\begin{enumerate}
\item[(a)] $H^i(Y, E_n(s))=0$ for any $i \in \Z_{>0}$, $n \in \Z_{>0}$, and $s \in \Z_{\geq s_2}$. 
\item[(b)] The induced homomorphism $R:E_{n+1}(s) \to E_n(s)$ is surjective 
for any $n \in \Z_{>0}$, and $s \in \Z_{\geq s_2}$. 
\end{enumerate}
By (1) and Proposition \ref{p-div-coherent}, 
$E_n(s)$ is a coherent $W_n\MO_Y$-module. 
Thus, (a) holds. 
Let us prove (b). 
We have the following exact sequence of $W\MO_X$-modules (Definition \ref{d-gr-def})
\[
0 \to {\rm gr}^n W\Omega_X^r \to W_{n+1}\Omega_X^r \to W_n \Omega_X^r \to 0. 
\]
Applying $(-) \otimes_{W\MO_X} W\MO_X(spA)$ and $R^if_*(-)$, 
we obtain the following exact sequence: 
\[
E_{n+1}(s) \xrightarrow{R} E_n(s) \to R^1f_*({\rm gr}^n W\Omega_X^N \otimes_{W\MO_X} W\MO_X(spA))
\]
Therefore, in order to prove (b), it is enough to find $s_3 \in \Z_{>0}$ 
which satisfies the following condition (c). 
\begin{enumerate}
\item[(c)] The equation 
$R^1f_*({\rm gr}^n W\Omega_X^N \otimes_{W\MO_X} W\MO_X(spA))=0$ 
holds for any $n \in \Z_{>0}$ and $s \in \Z_{\geq s_3}$. 
\end{enumerate}
By (1) and Proposition \ref{p-BZ}, 
it follows from induction on $n$ that there exists a positive integer $s_3$ such that 
\[
R^if_*(((F^e_X)_*\Omega_X^r) \otimes_{W\MO_X} W\MO_X(sA) )=0
\]
\[
R^if_*(((F^e_X)_*B_n\Omega_X^r) \otimes_{W\MO_X} W\MO_X(sA) )=0
\]
\[
R^if_*(((F^e_X)_*Z_n\Omega_X^r) \otimes_{W\MO_X} W\MO_X(sA) )=0
\]
for any $i \in \Z_{>0}$, $n, r, e \in \Z_{\geq 0}$ and $s \in \Z_{\geq s_3}$. 
Thus, we get 
\[
R^if_*\left(
\frac{(F_X^{n+1})_*\Omega_X^r}{(F_X)_*(B_n\Omega_X^r)} \otimes_{W\MO_X} W\MO_X(sA)
\right)=0,
\]
\[
R^if_*\left(
\frac{(F_X^{n+1})_*\Omega_X^r}{(F_X)_*(Z_n\Omega_X^r)} \otimes_{W\MO_X} W\MO_X(sA)
\right)=0.
\]
By the exact sequence (Proposition \ref{p-DRW-proj-system}(2))
\[
0 \to \frac{(F_X^{n+1})_*\Omega_X^r}{(F_X)_*B_n\Omega_X^r} \to 
(F_X)_*({\rm gr}^n W\Omega_X^r) \to 
\frac{(F_X^{n+1})_*\Omega_X^{r-1}}{(F_X)_*Z_n\Omega_X^{r-1}} \to 0, 
\]
it holds that 
\[
R^if_*((F_X)_*({\rm gr}^n W\Omega_X^r) \otimes_{W\MO_X} W\MO_X(sA))=0. 
\]
Then the assertion (c) holds by the following computation: 
\begin{eqnarray*}
0&=&R^if_*((F_X)_*({\rm gr}^n W\Omega_X^r) \otimes_{W\MO_X} W\MO_X(sA))\\
&\simeq & R^if_*(F_X)_* ({\rm gr}^n W\Omega_X^r \otimes_{W\MO_X} W\MO_X(psA))\\
& \simeq & (F_Y)_*R^if_*({\rm gr}^n W\Omega_X^r \otimes_{W\MO_X} W\MO_X(psA)), 
\end{eqnarray*}
where the last isomorphism follows from the fact that 
both $(F_X)_*$ and $(F_Y)_*$ are exact functors. 
This completes the proof of (3).  

Let us show (4). 
Fix a positive integer $t$ such that $p^t \geq \max \{s_1, ps_2\}$. 
We first prove 
\begin{equation}\label{e1-DRW-KV}
R^if_*(W\Omega_X^N \otimes_{W\MO_X} W\MO_X(p^tA))=0
\end{equation}
for any $i>0$. 
We obtain the following isomorphisms in the derived category of $W\MO_Y$-modules: 
\begin{eqnarray*}
&& Rf_*(W\Omega_X^N \otimes_{W\MO_X} W\MO_X(p^tA))\\
&\simeq& Rf_*((\varprojlim_n  W_n\Omega_X^N) \otimes_{W\MO_X} W\MO_X(p^tA))\\
&\simeq& Rf_* \varprojlim_n  (W_n\Omega_X^N \otimes_{W\MO_X} W\MO_X(p^tA))\\
&\simeq& Rf_* R\varprojlim_n  (W_n\Omega_X^N \otimes_{W\MO_X} W\MO_X(p^tA))\\
&\simeq& R\varprojlim_n  Rf_* (W_n\Omega_X^N \otimes_{W\MO_X} W\MO_X(p^tA))\\
&\simeq& R\varprojlim_n  f_* (W_n\Omega_X^N \otimes_{W\MO_X} W\MO_X(p^tA))\\
&\simeq& \varprojlim_n  f_* (W_n\Omega_X^N \otimes_{W\MO_X} W\MO_X(p^tA)),\\
\end{eqnarray*}
where the third isomorphism holds by the Mittag--Leffler condition and 
the fifth and the last isomorphisms 
follow from (2) and (3), respectively. 
This completes the proof of the equation (\ref{e1-DRW-KV}).

Then the assertion (4) holds by the following calculation: 
\begin{eqnarray*}
&&R^if_*(W\Omega_X^N \otimes_{W\MO_X} W\MO_X(A)) \\
&\simeq & R^if_*((F_X^t)_*(W\Omega_X^N) \otimes_{W\MO_X} W\MO_X(A))\\
&\simeq & R^if_*((F_X^t)_*(W\Omega_X^N \otimes_{W\MO_X} (F_X^t)^*W\MO_X(A)))\\
&\simeq & (F_Y^t)_*R^if_*(W\Omega_X^N \otimes_{W\MO_X} W\MO_X(p^tA))\\
&=& 0, 
\end{eqnarray*}
where the first isomorphism follows from Theorem \ref{t-descent}, 
the second one holds by the projection formula (\cite[Ch. 0, (5.4.10)]{Gro64}),
the third one follows from the fact that $(F_X^t)_*$ 
and $(F_Y^t)_*$ are exact functors, 
and the last equality holds by (\ref{e1-DRW-KV}). 
\end{proof}

In the proof of Theorem \ref{t-Kod-van}, 
${\rm gr}^n W\Omega_X^N$ plays a crucial role. 
In Remark \ref{r-gr-desc}, we give a description of ${\rm gr}^n W\Omega_X^N$ 
via the Grothendieck duality. 
The results in Remark \ref{r-gr-desc} will not be used in the rest of this paper. 

\begin{rem}\label{r-gr-desc}
Let $k$ be a perfect field of characteristic $p>0$ and 
let $X$ be an $N$-dimensional smooth variety over $k$. 
\begin{enumerate}
\item 
For any $r \in \Z_{\geq 0}$, 
we have the decomposition of $p : W_{n+1}\Omega^r_X \to W_{n+1}\Omega^r_X$ into the surjection $R$ and the injection $\underline p$ \cite[Proposition 3.4]{Ill79}: 
\begin{equation}\label{e1-r-gr-desc}
p : W_{n+1}\Omega^r_X \xrightarrow{R} W_n\Omega^r_X \xrightarrow{\underline p} W_{n+1}\Omega^r_X. 
\end{equation}
In particular, we obtain another exact sequence: 
\begin{equation}\label{e2-r-gr-desc}
0 \to W_n\Omega^r_X \xrightarrow{\underline p} W_{n+1}\Omega^r_X \to 
W_{n+1}\Omega^r_X/(p) \to 0. 
\end{equation}
\item 
We now show that $W_{n+1}\MO_X/(p)$ is a maximal Cohen--Macaulay $W_{n+1}\MO_X$-module. 
Fix a closed point $x \in X$. 
For $i < N$, we have that 
$H^i_x(W_n\MO_X) =0$ (cf. Theorem \ref{t-div-CM}). 
It is well known that the Frobenius homomorphism $F_X: \MO_{X, x} \to (F_X)_*\MO_{X, x}$ splits 
as an $\MO_X$-module homomorphism. 
Then we can show by induction on $n$ that 
\[
F_X : H^N_x(W_n\MO_X) \to H^N_x((F_X)_*W_n\MO_X)
\]
is injective. It follows from $\underline p = V \circ F_X$ that 
also $\underline p: H^N_x(W_n\MO_X) \to H^N_x(W_{n+1}\MO_X)$ is injective. 
By (\ref{e2-r-gr-desc}), we obtain 
$H^i_x(W_{n+1}\MO_X/(p))=0$ for any $i <N$, i.e. $W_{n+1}\MO_X/(p)$ is maximal Cohen--Macaulay. 
\item 
Applying $\mathcal Hom_{W_{n+1}\MO_X}(-, W_{n+1}\Omega_X^N)$ to 
\[
p : W_{n+1}\MO_X \xrightarrow{R} W_n\MO_X \xrightarrow{\underline p} W_{n+1}\MO_X, 
\]
we obtain the decomposition of $p : W_{n+1}\Omega^N_X \to W_{n+1}\Omega^N_X$ into the surjection $\underline{p}^*$ and the injection $R^*$
\[
p : W_{n+1}\Omega^N_X \xrightarrow{\underline{p}^*} W_n\Omega^N_X \xrightarrow{R^*} W_{n+1}\Omega^N_X. 
\]
Indeed, the injectivity of the dual $R^*$ of $R$ is automatic, 
whilst the surjectivity of the dual $\underline{p}^*$ of $\underline p$ holds by 
(2) and (\ref{e2-r-gr-desc}). 
The decomposition of $p$ into the surjection and the injection is unique 
up to isomorphisms, 
it follows from 
(\ref{e1-r-gr-desc}) that $\underline{p}^* = \sigma \circ R$ and $R^* =\underline p \circ \sigma^{-1}$ 
for some $W_n\MO_X$-linear automorphism $\sigma: W_n\Omega_X^N \to W_n\Omega_X^N$. 
\item 
By (2) and (\ref{e2-r-gr-desc}), 
we have that 
\[
0 \to \mathcal Hom_{W_{n+1}\MO_X}(W_{n+1}\MO_X/(p), W_{n+1}\Omega_X^N) \to 
W_{n+1}\Omega_X^N \xrightarrow{\underline{p}^* = \sigma \circ R} W_n\Omega_X^N \to 0. 
\]
By ${\rm gr}^n W\Omega_X^N = \Ker( R: W_{n+1}\Omega_X^N \to W_n\Omega_X^N)$ 
(Definition \ref{d-gr-def}), we obtain an isomorphism of $W_{n+1}\MO_X$-modules 
\[
{\rm gr}^n W\Omega_X^N \simeq  \mathcal Hom_{W_{n+1}\MO_X}(W_{n+1}\MO_X/(p), W_{n+1}\Omega_X^N). 
\]
\end{enumerate}
\end{rem}

\subsection{Splitting criteria}\label{ss4-split}

The purpose of this subsection is to prove that 
\[
W\MO_Y \to f_*W\MO_X 
\]
splits if a finite surjective morphism $f:X \to Y$ of smooth varieties 
satisfies one of the following properties. 
\begin{enumerate}
\item 
$K(X)/K(Y)$ is a Galois extension such that its extension degree $[K(X):K(Y)]$ is not divisible by $p$ (Lemma \ref{l-finite-trace}). 
\item $Y=X \times_k k'$ and $f: Y =X \times_k k' \to X$ is its projection, 
where $k \subset k'$ is a field extension of finite degree (Proposition \ref{p-pull-field-ext}). 
\end{enumerate}

\subsubsection{Galois covers with $|G|$ not divisible by $p$}

\begin{lem}\label{l-finite-trace}
Let $f:X \to Y$ be a finite surjective morphism of 
integral normal excellent $\F_p$-schemes. 
Assume that the induced field extension $K(Y) \subset K(X)$ is 
Galois and that the cardinality $|G|$ of its Galois group $G$ is not divisible by $p$. 
Let $D$ be an $\R$-divisor on $Y$. 
Then the following hold. 
\begin{enumerate}
\item 
The induced $W_n\MO_Y$-module homomorphism
$$W_n\MO_Y(D) \to f_*(W_n\MO_X(f^*D))$$ 
that is defined in Subsection \ref{ss-def-pullback-div} splits for any $n \in \Z_{>0}$. 
\item 
The induced $W\MO_Y$-module homomorphism  
$$W\MO_Y(D) \to f_*(W\MO_X(f^*D))$$ 
that is defined in Subsection \ref{ss-def-pullback-div} splits. 
\end{enumerate}
\end{lem}

\begin{proof}
We omit the proof of (1), as it is the same as the one of (2). 
Let us prove (2). 
Removing closed subsets of $X$ and $Y$ whose codimensions are at least two, 
the problem is reduced to the case when 
\begin{itemize}
\item both $X$ and $Y$ are regular, and 
\item if $D_1$ and $D_2$ are distinct prime divisors contained in $\Supp\,D$, 
then $D_1 \cap D_2=\emptyset$. 
\end{itemize}
For any $\sigma \in G$ and $n \in \Z_{\geq 0}$, 
the $K(Y)$-algebra automorphism $\sigma:K(X) \to K(X)$ 
induces an $\MO_Y$-module automorphism 
$$\sigma:f_*(\MO_X(p^nf^*D)) \to f_*\sigma_*(\MO_X(p^n\sigma^*f^*D))
=f_*(\MO_Y(p^nf^*D)).$$
Hence, the $W(K(Y))$-algebra automorphism $W(\sigma):W(K(X)) \to W(K(X))$ 
induces a $W\MO_Y$-module automorphism 
$$W(\sigma):f_*(W\MO_X(f^*D)) \to f_*(W\MO_X(f^*D)).$$
We denote by $(f_*W\MO_X(f^*D))^G$ the $G$-invariant 
$W\MO_Y$-submodule of $f_*(W\MO_X(f^*D))$, 
i.e. for any open subset $V$ of $Y$, we set 
$$\Gamma(V, f_*(W\MO_X(f^*D))^G):=\bigcap_{\sigma \in G} 
\{\beta \in \Gamma(V, f_*(W\MO_X(f^*D)))\,|\, W(\sigma)(\beta)=\beta\}.$$
In particular, we obtain $W\MO_Y$-module homomorphisms: 
$$W\MO_Y(D) \xrightarrow{\rho'} (f_*W\MO_X(f^*D))^G \hookrightarrow f_*W\MO_X(f^*D).$$

\setcounter{step}{0}
\begin{step}\label{s1-finite-trace}
The induced ring homomorphism 
$$\theta:W(K(Y)) \to W(K(X))^G$$
is bijective. 
\end{step}

\begin{proof}[Proof of Step \ref{s1-finite-trace}]
Since $\theta$ is automatically injective, 
it suffices to prove that $\theta$ is surjective. 
Take 
$$b=(b_0, b_1, \cdots) \in W(K(X))^G,$$
i.e. $b \in W(K(X))$, $b_i \in K(X)$, and $W(\sigma)(b)=b$ for any $\sigma \in G$. 
Then we have that 
$$(b_0, b_1, \cdots)=W(\sigma)(b)=(\sigma(b_0), \sigma(b_1), \cdots),$$
hence $\sigma(b_i)=b_i$ for any $i\geq 0$ and $\sigma \in G$. 
Therefore, $b_i \in K(Y)$ for any $i \geq 0$. 
This completes the proof of  Step \ref{s1-finite-trace}. 
\end{proof}

\begin{step}\label{s2-finite-trace}
The induced $\MO_Y$-module homomorphism 
$$\rho:\MO_Y(D) \to (f_*\MO_X(f^*D))^G$$ 
is an isomorphism. 
\end{step}

\begin{proof}[Proof of Step \ref{s2-finite-trace}]
Replacing $D$ by a slightly larger $\Q$-divisor, 
we may assume that $D$ is a $\Q$-divisor (cf. Remark \ref{r-perturb}). 
It is obvious that $\rho$ is injective, 
hence let us prove that $\rho$ is surjective. 
Since the problem is local on $Y$, we may assume that 
$$D=q {\rm div}(\psi)$$
for some $q \in \Q \setminus \{0\}$ and $\psi \in K(Y)^{\times}$. 
Take 
$$\varphi \in \Gamma(X, f^*D)^G=\Gamma(X, f^*D) \cap K(Y).$$ 
We have that 
$$f^*({\rm div}(\varphi)+q{\rm div}(\psi)) \geq 0.$$
Pick $m \in \Z_{>0}$ such that $mq \in \Z$. 
It holds that 
$$\varphi^m\psi^{mq} \in \Gamma(X, \MO_X) \cap K(Y)=\Gamma(Y, \MO_Y).$$
This implies that 
$${\rm div}(\varphi)+q{\rm div}(\psi) \geq 0,$$
which is equivalent to 
$$\varphi \in \Gamma(X, D).$$
Thus, $\rho$ is surjective. 
This completes the proof of Step \ref{s2-finite-trace}. 
\end{proof}

\begin{step}\label{s3-finite-trace}
The induced $W\MO_Y$-module homomorphism: 
$$\rho':W\MO_Y(D) \to (f_*W\MO_X(f^*D))^G$$
is an isomorphism. 
\end{step}

\begin{proof}[Proof of Step \ref{s3-finite-trace}]
It is clear that $\rho'$ is injective, 
hence let us prove that $\rho'$ is surjective. 
Take $\beta \in \Gamma(Y, (f_*W\MO_X(f^*D))^G)$. 
We have that 
\[
\beta=(b_0, b_1, \cdots) \in W(K(X))^G=W(K(Y)), 
\]
where the last equality is guaranteed by Step \ref{s1-finite-trace}. 
Then, for any $n \in \Z_{\geq 0}$,  it holds that 
$$b_n \in K(Y) \cap \Gamma(X, \MO_X(p^nf^*(D)))=\Gamma(Y, \MO_X(p^n D)),$$
where the last equality follows from Step \ref{s2-finite-trace}. 
Therefore, $\rho'$ is surjective. 
This completes the proof of Step \ref{s3-finite-trace}. 
\end{proof}

\begin{step}\label{s4-finite-trace}
The assertion (2) holds. 
\end{step}

\begin{proof}[Proof of Step \ref{s4-finite-trace}]
We have an $W\MO_Y$-module homomorphism: 
\begin{eqnarray*}
T:f_*(W\MO_X(f^*D)) &\to & W\MO_Y(D), \\
\beta &\mapsto& \frac{1}{|G|}
\sum_{\sigma \in G} \sigma(\beta).
\end{eqnarray*}
Then $T$ gives the splitting of 
the natural $W\MO_X$-module homomorphism $W\MO_Y(D) \to f_*(W\MO_X(f^*D))$. 
Thus (2) holds, which completes the proof of Step \ref{s4-finite-trace}. 
\end{proof}
Step \ref{s4-finite-trace} completes the proof of Lemma \ref{l-finite-trace}. 
\end{proof}

\begin{rem}
With notation as in the statement of Lemma \ref{l-finite-trace}, 
the $W\MO_{Y, \Q}$-module homomorphism (for the definition of $W\MO_{Y, \Q}$, see Subsection \ref{ss-notation}(\ref{ss-notation-Q}))  
$$\gamma:W\MO_Y(D)_{\Q} \to f_*(W\MO_X(f^*(D)))_{\Q}$$
splits. 
Furthermore, this splitting holds even if 
we drop the assumption that $f$ is Galois. 
Although we do not use this fact in the paper, 
let us give a sketch of a proof. 
If $f$ is purely inseparable, then $\gamma$ is an isomorphism automatically. 
If $f$ is Galois, then the same proof as above works 
even for the case when $|G|$ is divisible by $p$. 
If $f$ is separable, then the problem is reduced, by taking a splitting field, 
to the case when $f$ is Galois. 
\end{rem}

\subsubsection{Extension of the base field}

\begin{lem}\label{l-div-etale}
Let $f:X \to Y$ be an \'etale morphism of regular excellent $\F_p$-schemes. 
Let $D$ be an $\R$-divisor on $Y$. 
Then the induced $W_n\MO_X$-module homomorphism 
\[
(W_nf)^*(W_n\MO_Y(D)) \to W_n\MO_X(f^*D)
\]
is an isomorphism, 
where $W_nf$ denotes the induced morphism $W_nX \to W_nY$ of schemes. 
\end{lem}

\begin{proof}
Recall that the $W_n\MO_X$-module homomorphism 
$$W_n\MO_Y(D) \to f_*W_n\MO_X(f^*D)=(W_nf)_*(W_n\MO_X(f^*D))$$ 
defined in Subsection \ref{ss-def-pullback-div} induces, by adjunction, a $W_n\MO_Y$-module homomorphism 
$$\alpha_n:(W_nf)^*(W_n\MO_Y(D)) \to W_n\MO_X(f^*D).$$
It follows from \cite[Ch. 0, Proposition 1.5.8]{Ill79} that 
$$(F_X^e)_*(W_nf)^*(W_n\MO_Y(D))\simeq (W_nf)^*((F_Y^e)_*W_n\MO_Y(D)).$$
Thus we get a commutative diagram:  
{\tiny
\[
\begin{CD}
0 @>>> (F^{e+1}_X)_*(W_nf)^*(W_n\MO_Y(pD)) @>>> (F_X^e)_*(W_{n+1}f)^*(W_{n+1}\MO_Y(D)) @>>> (F_X^e)_*f^*(\MO_Y(D)) @>>> 0\\
@. @VV(F^{e+1}_X)_*\alpha_nV @VV(F^{e}_X)_*\alpha_{n+1}V @VV(F^{e}_X)_*\alpha_1 V\\
0 @>>> (F^{e+1}_X)_*W_n\MO_X(pf^*D)) @>>> (F_X^e)_*W_{n+1}\MO_X(f^*D) @>>> (F_X^e)_*\MO_X(f^*D) @>>> 0,\\
\end{CD}
\]
}
where both the horizontal sequences are exact (cf. Proposition \ref{p-div-induction}). 
By the snake lemma and induction on $n$, 
it suffices to prove that 
$$\alpha_1:f^*\MO_Y(D) \to \MO_X(f^*D).$$ 
is an isomorphism. 
Since $f$ is \'etale, $f^*S$ is a reduced divisor for any prime divisor $S$ of $Y$. 
Thus we have that 
$$\MO_Y(D)=\MO_Y(\llcorner D \lrcorner) \quad {\rm and} \quad 
\MO_X(f^*D)=\MO_X(\llcorner f^*D\lrcorner)=\MO_X(f^*(\llcorner D\lrcorner)).$$
Replacing $D$ by $\llcorner D \lrcorner$, we may assume that $D$ is a Cartier divisor. 
In this case, we can check directly that $\alpha_1$ is an isomorphism. 
\end{proof}

\begin{prop}\label{p-pull-field-ext}
Let $k \subset k'$ be a finite extension of perfect fields of characteristic $p>0$. 
Let $X$ be a geometrically connected smooth variety over $k$. 
Set $X':=X \times_k k'$ and let $f:X' \to X$ be the induced morphism. 
Let $D$ be an $\R$-divisor on $X$ and set $D':=f^*D$. 
Fix $e \in \Z_{\geq 0}$. 
Then the following hold. 
\begin{enumerate}
\item 
The induced $W_n\MO_X$-module homomorphism 
\[
\gamma_{e, n}:
((F^e_X)_*W_n\MO_{X}(D)) \otimes_{W(k)} W(k') \to 
(F^e_X)_*f_*(W_n\MO_{X'}(D'))
\]
is an isomorphism for any $n \in \Z_{>0}$. 
\item 
The induced $W\MO_X$-module homomorphism 
\[
\gamma_e: ((F^e_X)_*W\MO_{X}(D)) \otimes_{W(k)} W(k') \to 
(F^e_X)_*f_*(W\MO_{X'}(D'))
\]
is an isomorphism. 
\item 
The induced $W\MO_X$-module homomorphism 
\[
W\MO_{X}(D) \to f_*(W\MO_{X'}(D'))
\]
splits. 
\end{enumerate}
\end{prop}

\begin{proof}
Let us show (1). 
It follows from \cite[Ch. 0, Proposition 1.5.8]{Ill79} that 
\[
(F_{X'}^e)_*(W_nf)^*(W_n\MO_X(D))\simeq (W_nf)^*((F_X^e)_*W_n\MO_X(D)).
\]
Then the assertion (1) follows from Lemma \ref{l-div-etale}.

Let us show (2). 
By (1), we obtain an isomorphism: 
\[
\varprojlim_n 
\gamma_{e, n}:
\varprojlim_n (((F^e_X)_*W_n\MO_{X}(D)) \otimes_{W(k)} W(k')) \to 
\varprojlim_n (F^e_X)_*f_*(W_n\MO_{X'}(D')). 
\]
It holds that 
$\varprojlim_n (F^e_X)_*f_*(W_n\MO_{X'}(D')) \simeq (F^e_X)_*f_*(W\MO_{X'}(D'))$. 
Since $W(k')$ is a free $W(k)$-module of finite rank by \cite[Theorem 8.4]{Mat89}, 
we have that $\varprojlim_n$ and $(-) \otimes_{W(k)} W(k')$ commutes: 
\[
\varprojlim_n (((F^e_X)_*W_n\MO_{X}(D)) \otimes_{W(k)} W(k'))
\simeq 
((F^e_X)_*W\MO_{X}(D)) \otimes_{W(k)} W(k'). 
\]
Thus, (2) holds. The assertion (3) follows directly from (2). 
\end{proof}

\subsection{$\mathcal Hom(f_*W\MO_Y(D), -)$ vs $\mathcal Hom(f_*W_n\MO_Y(D), -)$}\label{ss-WO-WO}

The purpose of this subsection is to prove Proposition \ref{p-WO-vs-WO-Y}, 
which assures 
that the natural $W\MO_X$-module homomorphism 
\[
\mathcal Hom_{W_n\MO_X} (f_*W_n\MO_Y(f^*D), W_n\Omega_X^N) \to 
\mathcal Hom_{W\MO_X}(f_*W\MO_Y(f^*D), W_n\Omega_X^N) 
\]
is an isomorphism under the assumption that $f$ is weakly $\ell$-cyclic in the sense of Definition \ref{d-weak-l-cyclic}. 
The key result is Lemma \ref{l-WO-vs-WO-Y}.

\begin{dfn}\label{d-weak-l-cyclic}
Let $\ell$ be a positive integer. 
\begin{enumerate}
\item 
We say that a ring extension $A \subset B$ is {\em weakly} $\ell$-{\em cyclic} if 
there exist $y_1, ..., y_r \in B$ such that 
\begin{enumerate}
\item $B=A y_1+ \cdots +A y_r$, and 
\item $y_1^{\ell} \in A, ..., y_r^{\ell} \in A$. 
\end{enumerate}
In particular, $B$ is a finitely generated $A$-module. 
\item 
Let $f:Y \to X$ be a finite surjective morphism of reduced schemes. 
We say that $f$ is 
{\em weakly} $\ell$-{\em cyclic} if there exists an affine open cover $X=\bigcup_{i \in I} X_i$ such that, for any $i \in I$, the induced ring extension $\MO_X(X_i) \hookrightarrow \MO_Y(f^{-1}(X_i))$ 
is weakly $\ell$-cyclic in the sense of (1) (cf. Lemma \ref{l-weak-l-cyclic2}). 
\end{enumerate}
\end{dfn}

\begin{lem}\label{l-weak-l-cyclic2}
Let $\ell$ be a positive integer and let $f:Y \to X$ be a finite surjective morphism of reduced schemes. 
Then the following are equivalent. 
\begin{enumerate}
\item $f$ is weakly $\ell$-cyclic. 
\item For an arbitrary affine open cover $X=\bigcup_{i \in I} X_i$, 
the induced ring extension $\MO_X(X_i) \hookrightarrow \MO_Y(f^{-1}(X_i))$ is weakly $\ell$-cyclic in the sense of Definition \ref{d-weak-l-cyclic}(1). 
\end{enumerate}
\end{lem}

\begin{proof}
It suffices to show that (1) implies (2). 
After replacing $X$ and an open cover suitably, it is enough to prove $(*)$. 
\begin{enumerate}
\item[$(*)$] Assume that $X = \Spec\,A$ and $Y = \Spec\,B$. 
Pick elements $f_1, ..., f_n \in A$ such that $Af_1 + \cdots + Af_n =A$.  
If $A_{f_i} \hookrightarrow B_{f_i}$ is weakly $\ell$-cyclic, 
then also $A \hookrightarrow B$ is weakly $\ell$-cyclic. 
\end{enumerate}

Fix  $1 \leq i \leq n$. 
It follows from Definition \ref{d-weak-l-cyclic}(1) that 
there exist $y_{i, 1}, ...., y_{i, r_i} \in B_{f_i}$ such that 
\begin{enumerate}
\item[(a$)_i$] $B_{f_i}=A_{f_i} y_{i, 1}+ \cdots +A_{f_i} y_{i, r_i}$, and 
\item[(b$)_i$] $y_{i, 1}^{\ell} \in A_{f_i}, ..., y_{i, r_i}^{\ell} \in A_{f_i}$. 
\end{enumerate}
After replacing $y_{i, j}$ by $f_i^{\ell N}y_{i, j}$ for suitable $N \in \Z_{>0}$, 
we may assume that 
\begin{enumerate}
\item[(a$)'_i$] $f_i^{M}B \subset A y_{i, 1}+ \cdots +A y_{i, r_i}$ for some $M \in \Z_{>0}$, 
\item[(b$)'_i$] $y_{i, 1}^{\ell} \in A, ..., y_{i, r_i}^{\ell} \in A$, and 
\item[(c$)'_i$] $y_{i, 1}, ...., y_{i, r_i}  \in B$. 
\end{enumerate}
We have $a_1 f^M_1 + \cdots + a_nf^M_n =1$ for some $a_1, ..., a_n \in A$. 
By (a$)'_i$, it holds that $B = \sum_{i, j} A y_{i, j}$. 
This equation and  (b$)'_i$ imply that $A \hookrightarrow B$ is weakly $\ell$-cyclic. 
\end{proof}

\begin{lem}\label{l-WO-vs-WO-Y}
Fix $\ell \in \Z_{>0} \setminus p\Z$ and $n \in \Z_{>0}$. 
Let $f:Y \to X$ be a finite surjective morphism of 
affine excellent integral normal $\F_p$-schemes. 
Set $A:=\MO_X(X)$ and $B:=\MO_Y(Y)$. 
Assume that the induced ring extension $A \hookrightarrow B$ is weakly $\ell$-cyclic. 
Let $D$ be an $\R$-divisor on $X$. 
Take $\zeta \in \Gamma(Y, V^n(W\MO_Y(f^*D)))$. 
Then there exists a nonzero element $a \in A \setminus \{0\}$ such that 
\[
\underline{a} \zeta \in \Gamma(X, V^n(W\MO_X)) \cdot \Gamma(Y, W\MO_Y(f^*D)). 
\]
\end{lem}

\begin{proof}
We first reduce the problem to the case when $D=0$. 
Since $Y$ is affine, there exists a nonzero element $a_1 \in A \setminus \{0\}$ 
such that $\underline{a_1} \zeta \in \Gamma(Y, V^n(W\MO_Y))$. 
Since we assume that the assertion holds when $D=0$, 
there exists $a_2 \in A \setminus \{0\}$ such that 
\[
\underline{a_1a_2} \zeta \in \Gamma(X, V^n(W\MO_X)) \cdot \Gamma(Y, W\MO_Y).
\]
Since $Y$ is affine, there exists a nonzero element $a_3 \in A \setminus \{0\}$ 
such that 
\[
\underline{a_1a_2a_3} \zeta \in \Gamma(X, V^n(W\MO_X)) \cdot \Gamma(Y, W\MO_Y(f^*D)).
\]
Set $a:=a_1a_2a_3$. 
Then we completes the proof of the reduction to the case when $D=0$. 

From now on, we treat the case when $D=0$. 
If $\ell=1$, then there is nothing to show. 
We assume that $\ell \geq 2$. 

Fix $m \in \Z_{> 0}$.  
We now show that there exist two integers $u_m, v_m \in \Z$ which satisfy 
the following properties (i)--(iii) 
(such $u_m$ and $v_m$ are unique, although we do not use this fact). 
\begin{enumerate}
\renewcommand{\labelenumi}{(\roman{enumi})}
\item $0< u_m  < p^m$. 
\item $0< v_m < \ell$. 
\item $-u_m\ell + v_m p^m=1$. 
\end{enumerate}
By $\ell \in \Z_{\geq 2} \setminus p \Z$, there 
exists $v_m \in \Z$ such that 
$0 < v_m < \ell$ and $v_m p^m \equiv 1 \mod \ell$. 
In particular, (ii) holds. 
We define $u_m \in \Q$ by the equation (iii). 
By $v_m p^m \equiv 1 \mod \ell$, we obtain $u_m \in \Z$. 
Furthermore, we have $u_m \ell =v_mp^m -1$. 
By $v_m \in \Z_{>0}$, we get $u_m \ell =v_mp^m -1 > 0$, which implies $u_m >0$. 
Then we obtain the remaining property (i) by the following computation: 
\[
u_m = \frac{1}{\ell}\left(v_mp^m -1 \right) < \frac{1}{\ell}\left(\ell p^m -1 \right) < p^m,
\]
where the equality follows from (iii) and the first inequality holds by (ii). 
This completes the proof of (i)--(iii). 

Since $A \subset B$ is weakly $\ell$-cyclic, 
there exist $y_1, ..., y_r \in B$ such that 
\begin{enumerate}
\item[(a)] $B=Ay_1+\cdots + A y_r$, and 
\item[(b)] $y_1^{\ell}, ..., y_r^{\ell} \in A$. 
\end{enumerate}
By (a) and $\zeta \in V^n(WB)$, we can write 
\[
\zeta = \sum_{j=1}^r V^n(\underline{a_{1j} y_j}) + \sum_{j=1}^r V^{n+1}(\underline{a_{2j} y_j}) + 
\cdots 
=\sum_{m=n}^{\infty} \sum_{j=1}^r V^m(\underline{a_{mj} y_j})
\]
for some $a_{mj} \in A$. 
Set $x_j:=y_j^{\ell} \in A$ (cf. (b)). 
Then we have that 
\[
y_j x_j^{u_m} = y_j \cdot y_j^{u_m\ell} = y_j^{1+u_m\ell} =  y_j^{v_mp^m}. 
\]
By $p^m > u_m$, 
we get $y_jx_j^{p^m} = y_j x_j^{u_m} x_j^{p^m-u_m}=y_j^{v_mp^m} x_j^{p^m-u_m}$. 
Using this equation, it holds that 
\[
\underline{x_j}V^m(\underline{a_{mj} y_j}) 
= V^m(\underline{a_{mj} y_jx_j^{p^m}}) = 
V^m(\underline{a_{mj} y_j^{v_mp^m} x_j^{p^m-u_m}}) =
 \underline{y_j^{v_m}} V^m(\underline{a_{mj} x_j^{p^m-u_m}}). 
\]
For 
\[
z_j :=\frac{x_1 \cdots x_r}{x_j} = x_1 \cdots x_{j-1} x_{j+1} \cdots x_r \in A,
\]
we have that 
\begin{eqnarray*}
\underline{x_1 \cdots x_r} \cdot \zeta 
&=& \underline{x_1 \cdots x_r}  \sum_{m=n}^{\infty} \sum_{j=1}^r V^m(\underline{a_{mj} y_j})\\
&=&  \sum_{j=1}^r   \sum_{m=n}^{\infty}\underline{x_1 \cdots x_r}V^m(\underline{a_{mj} y_j})\\
&=&  \sum_{j=1}^r  \sum_{m=n}^{\infty} \underline{z_jx_j} V^m(\underline{a_{mj} y_j})\\
&=&  \sum_{j=1}^r \sum_{m=n}^{\infty}\underline{z_jy_j^{v_m}} V^m(\underline{a_{mj} x_j^{p^m-u_m}}).
\end{eqnarray*}
Recall that $v_m \in \{1, 2,..., \ell-1\}$. 
For each $v \in \{1, 2,..., \ell-1\}$, we set 
$M_v := \{m \in \Z_{\geq n}\,|\, v_m = v\}$. 
Then, for $a := x_1 \cdots x_r$,  we obtain 
\begin{eqnarray*}
\underline{a} \cdot \zeta &=& \underline{x_1 \cdots x_r} \cdot \zeta 
=\sum_{j=1}^r \sum_{v=1}^{\ell-1} \xi_{vj}\\
\xi_{vj} &:=& \sum_{m \in M_v} \underline{z_jy_j^v} V^m(\underline{a_{mj} x_j^{p^m-u_m}})\\
&=& \underline{z_jy_j^v} \sum_{m \in M_v} V^m(\underline{a_{mj} x_j^{p^m-u_m}}) \in (WB) \cdot V^n(WA),
\end{eqnarray*}
as desired. 
\end{proof}

\begin{prop}\label{p-WO-vs-WO-Y}
Fix $\ell \in \Z_{>0} \setminus p\Z$ and $n \in \Z_{>0}$. 
Let $f:Y \to X$ be a finite surjective morphism of excellent integral normal $\F_p$-schemes. 
Assume that $f$ is weakly $\ell$-cyclic. 
Let $D$ be an $\R$-divisor on $X$. 
Let $M$ be a coherent $W_n\MO_X$-module such that 
the induced map $M(U) \to M_{\xi}$ is injective 
for any non-empty open subset $U$ of $X$, 
where $M_{\xi}$ denotes the stalk at the generic point $\xi$ of $X$. 
Then the induced $W\MO_X$-module homomorphism 
\[
\theta:\mathcal Hom_{W_n\MO_X} (f_*W_n\MO_Y(f^*D), M) \to 
\mathcal Hom_{W\MO_X}(f_*W\MO_Y(f^*D), M) 
\]
is an isomorphism. 
\end{prop}

\begin{proof}
It is clear that $\theta$ is injective. 
It suffices to show that $\theta$ is surjective. 
Since the problem is local, we may assume that 
$X$ and $Y$ are affine and that 
the ring extension $\MO_X(X) \hookrightarrow \MO_Y(Y)$ is weakly $\ell$-cyclic. 
Take $\psi \in \Gamma(X, \mathcal Hom_{W\MO_X} (f_*W\MO_Y(f^*D), M))$, i.e. 
$\psi:f_*W\MO_Y(f^*D)\to M$ is a $W\MO_X$-module homomorphism. 
It is enough to show that $\psi(f_*V^n(W\MO_Y(f^*D)))=0$. 
Take 
\[
\zeta \in \Gamma(X, f_*V^n(W\MO_Y(f^*D)))=\Gamma(Y, V^n(W\MO_Y(f^*D))).
\] 
By Lemma \ref{l-WO-vs-WO-Y}, there exists $a \in \MO_X(X) \setminus \{0\}$ such that 
\[
\underline{a} \zeta \in \Gamma(X, V^n(W\MO_X)) \cdot \Gamma(Y, W\MO_Y(f^*D)). 
\]
Therefore, we get $\psi(\underline{a} \zeta)=0$ in $M(X)$. 
Since $M(X) \to M_{\xi}$ is inietive, we get $\psi(\zeta)=0$, as desired. 
\end{proof}

\subsection{Kawamata--Viehweg vanishing}\label{ss-Kaw-Vie}

The main purpose of this subsection is to prove Theorem \ref{t-KVV}. 
We first recall a slightly modified version of Kawamata's covering trick 
(Lemma \ref{l-Kawamata-cover}, Theorem \ref{t-Kawamata-cover}, Proposition \ref{p-Kawamata}). 
We then prove Theorem \ref{t-KVV}. 

Although the proofs of Lemma \ref{l-Kawamata-cover} and Theorem \ref{t-Kawamata-cover} 
are very similar to the ones of \cite[Lemma 1-1-2]{KMM87} and \cite[Theorem 1-1-1]{KMM87} respectively, 
we give proofs for the sake of completeness. 

\begin{lem}\label{l-Kawamata-cover}
Let $k$ be a perfect field of characteristic $p>0$. 
Let $A$ be a $d$-dimensional regular local $k$-algebra such that 
the induced field extension $k \hookrightarrow A/\m$ is of finite degree, 
where $\m$ denotes the maximal ideal of $A$.  
Let $\{z_1, ..., z_d\}$ be a regular system of parameter of $A$. 
Fix $e \in \Z$ with $0 \leq e \leq d$. 
Let $u_1, ..., u_s$ be units of $A$. 
Fix $\ell \in \Z_{>0} \setminus p\Z$. 
Then the ring 
\[
B := A[x_1, ..., x_e, y_1, ..., y_s]/(x_i^{\ell} - z_i, y_j^{\ell} - u_j)_{i, j} 
\]
is a regular ring. 
Furthermore, if $\m_B$ is a maximal ideal of $B$, 
then $\{ \overline{x}_1, ..., \overline{x}_e, z_{e+1}, ..., z_d\}$ 
is a regular system of parameter of $B_{\m_B}$, 
where $\overline{x}_1, ..., \overline{x}_e \in B$ 
denote the images of $x_1, ..., x_e$, respectively. 
\end{lem}

\begin{proof}
We first reduce the problem to the case when $k \to A/\m$ is an isomorphism and $k$ has 
a primitive $\ell$-th root of the unity. 
Fix a maximal ideal $\m_B$ of $B$. 
Take a field extension $A/\m \hookrightarrow k'$ of finite degree 
such that 
\begin{itemize}
\item 
$k'$ contains  
a primitive $\ell$-th root of the unity, and 
\item 
there is a maximal ideal $\mathfrak n$ of $A \otimes_k k'$ such that 
the induced ring homomorphism $k' \to (A \otimes_k k') / \mathfrak n$ is an isomorphism. 
\end{itemize}
Set $A' := (A \otimes_k k')_{\mathfrak n}$ and $\m' := \mathfrak n A'$. 
Since $A \to A'$ is an essentially \'etale local homomorphism of local rings, 
$A \to A'$ is faithfully flat and $\{z_1, ..., z_d\}$ is a regular system of parameter of $A'$. 
In particular, also 
\[
B \to B' := B \otimes_A A' = A'[x_1, ..., x_e, y_1, ..., y_s]/(x_i^{\ell} - z_i, y_j^{\ell} - u_j)_{i, j}
\]
is faithfully flat and essentially \'etale. 
Therefore, there exists a maximal ideal $\m_{B'}$ lying over $\m_B$. 
As we are assuming that the assertion holds for $(A', \m', B')$, 
we have that 
\[
\m_{B'}B'_{\m_{B'}} = 
(\overline{x}_1, ... , \overline{x}_e, z_{e+1}, ...,  z_d) B'_{\m_{B'}}.
\]
Since $B_{\m_B} \to B'_{\m_{B'}}$ is an essentially \'etale local homomorphism, 
it holds that 
\[
\m_{B}B_{\m_{B}} = (\overline{x}_1, ... , \overline{x}_e, z_{e+1}, ...,  z_d) B_{\m_{B}}, 
\]
as required. 


It is enough to prove the assertion under the condition that 
$k \to A/\m$ is an isomorphism
and 
$k$ has a primitive $\ell$-th root of the unity $\zeta$. 
We have $u_1 - \alpha_1, ..., u_s -\alpha_s \in \m$ 
for some $\alpha_1, ..., \alpha_s \in k^{\times}$. 
Let $\overline{y}_1, ..., \overline y_s \in B$ be the images of $y_1, ..., y_s$, respectively. 
It is clear that $B$ is a finitely generated $A$-module. 
Hence, the pullback of $\m_B$ to $A$ is equal to $\m$. 
We obtain $k \xrightarrow{\simeq} A/\m \xrightarrow{\simeq} B/\m_B$. 
Then we can find $\beta_1, ..., \beta_s \in k$ such that $\overline y_j - \beta_j \in \m_B$ for any $1 \leq j \leq s$. 
We then have 
\[
\m_B =( \overline{x}_1, ..., \overline{x}_e, z_{e+1}, ..., z_d , 
\overline{y}_1 -\beta_1, ..., \overline{y}_s - \beta_s). 
\]
For any $1 \leq j \leq s$, 
we obtain $\beta_j^{\ell} = \overline{y}_j^{\ell} = u_j = \alpha_j$ in $B/\m_B$, 
hence we have $\beta_j^{\ell} = \alpha_j$ in $k$. 
By $\alpha_j \neq 0$, we have $\beta_j \neq 0$.

Fix $1 \leq j \leq s$. 
In order to show that 
$\{ \overline{x}_1, ..., \overline{x}_e, z_{e+1}, ..., z_d\}$ is a regular system of parameter of $B_{\m_B}$, it suffices to prove 
\[
\overline y_j - \beta_j \in ( \overline{x}_1, ..., \overline{x}_e, z_{e+1}, ..., z_d) B_{\m_B}. 
\]
For any $1 \leq m \leq \ell-1$, we have 
\[
\overline{y}_j -\zeta^m \beta_j =  (\overline y_j -\beta_j) + \beta_j ( 1-\zeta^m) 
\in \m_B + k^{\times} \subset (B_{\m_B})^{\times}. 
\]
Hence, we obtain 
\[
\overline{y}_j - \beta_j 
= \frac{\prod_{m=0}^{\ell} (\overline{y}_j -\zeta^m \beta_j)}{\prod_{m=1}^{\ell} (\overline{y}_j -\zeta^m \beta_j) } 
= \frac{\overline{y}_j^{\ell} - \alpha_j}{\prod_{m=1}^{\ell} (\overline{y}_j -\zeta^m \beta_j) } = 
\frac{u_j - \alpha_j}{\prod_{m=1}^{\ell} (\overline{y}_j -\zeta^m \beta_j) } 
\]
\[
\in \m B_{\m_B}
= (z_1, ..., z_e, z_{e+1}, ..., z_d) B_{\m_B} \subset 
( \overline{x}_1, ..., \overline{x}_e, z_{e+1}, ..., z_d) B_{\m_B}, 
\]
as required. 
\end{proof}

\begin{thm}\label{t-Kawamata-cover}
Fix $\ell \in \Z_{>0} \setminus p\Z$. 
Let $k$ be an infinite perfect field of characteristic $p>0$ such that 
$k$ has a primitive $\ell$-th root of the unity. 
Let $X$ be a quasi-projective smooth variety over $k$. 
Let $D$ be a $\Q$-divisor such that 
$\{D\}$ is simple normal crossing and $\ell D$ is a $\Z$-divisor. 
Then there exists a finite surjective $k$-morphism 
\[
g:Y \to X
\]
of smooth $k$-varieties that satisfy the following properties. 
\begin{enumerate}
\item 
$g$ is weakly $\ell$-cyclic. 
\item 
The induced field extension $K(Y)/K(X)$ is a Galois extension 
such that $[K(Y):K(X)]$ is not divisible by $p$. 
\item  $g^*D$ is a $\Z$-divisor. 
\end{enumerate}
\end{thm}

\begin{proof}
Replacing $k$ by its algebraic closure in $K(X)$, 
we may assume that $X$ is geometrically integral over $k$ 
\cite[(4.3.1), Proposition 4.5.9, Corollaire 4.6.11]{Gro65}. 
Set $d:=\dim X$. 
Let $\{ D \}= \sum_{i \in I} a_i\Gamma_i$ be the prime decomposition. 
Let $M$ be a very ample Cartier divisor on $X$ such that $\ell M - \Gamma_i$ is very ample for any $i \in I$. 
Take general members $H_j^{(i)} \in |\ell M - \Gamma_i|$ for $i \in I$ and $1 \leq j \leq d$. 
By the Bertini theorem \cite[Corollaire 6.7]{Jou83}, 
we may assume that 
\begin{enumerate}
\item[(i)] 
each $H_j^{(i)}$ is a prime divisor, and 
\item[(ii)]
$\{D\}+\sum_{i,j} H_j^{(i)}$ is a simple normal crossing divisor with $H_j^{(i)} \not\subset \Supp\,\{D\}$. 
\end{enumerate}
Pick an affine open cover $X=\bigcup_{\alpha \in A} U_{\alpha}$ 
with the transition functions 
$a_{\alpha \beta} \in H^0(U_{\alpha} \cap U_{\beta}, \MO_X^{\times})$ of $M$ and 
local sections $\varphi^{(i)}_{j, \alpha} \in H^0(U_{\alpha}, \MO_X)$ 
such that $(H^{(i)}_j + \Gamma_i)|_{U_{\alpha}} = {\rm div}(\varphi^{(i)}_{j, \alpha})$ 
and that $\varphi^{(i)}_{j, \alpha} = a_{\alpha \beta}^{\ell} \varphi^{(i)}_{j, \beta}$. 
We have $K(X)[(\varphi^{(i)}_{j, \alpha})^{1/\ell}]_{i, j} = K(X)[(\varphi^{(i)}_{j, \beta})^{1/\ell}]_{i, j}$ for any $\alpha, \beta \in A$. 
Set 
\[
L:=K(X)[(\varphi^{(i)}_{j, \alpha})^{1/\ell}]_{i, j} 
\]
for some (hence any) $\alpha \in A$. 
Let $f:Y \to X$ be the normalisation of $X$ in $L$. 
Then (2) holds, because $L=K(Y)/K(X)$ is a Kummer extension such that $[K(Y):K(X)]$ is a divisor of $\ell^t$ for some $t \in \Z_{>0}$.

We now prove that the following property (iii) holds. 
\begin{enumerate}
\renewcommand{\labelenumi}{(\roman{enumi})}
\item[(iii)] $L = K(X)[t_{i, j}]_{i, j}/(t_{i, j}^{\ell} - \varphi^{(i)}_{j, \alpha})_{i, j}$. 
\end{enumerate}
Fix $i_1$ and $j_1$. 
Then $L_1:=K(X)[t_{i_1, j_1}]/(t^{\ell}_{i_1, j_1}-\varphi^{(i_1)}_{j_1, \alpha})$ is a field. 
Indeed, if $L_1$ is not a field, then $t^{\ell}_{i_1, j_1}-\varphi^{(i_1)}_{j_1, \alpha}$ 
is not irreducible over $K(X)$. 
Then there exist $\psi_{\alpha} \in K(X)$ and a divisor $\ell'$ of $\ell$ such that $1 < \ell' \leq \ell$ and 
$\varphi^{(i_1)}_{j_1, \alpha} = \psi^{\ell'}_{\alpha}$. 
By $\varphi^{(i)}_{j, \alpha} = a_{\alpha \beta}^{\ell} \varphi^{(i)}_{j, \beta}$, 
it holds that, for any $\beta \in A$, 
there exists $\psi_{\beta} \in K(X)$ such that 
$\varphi^{(i_1)}_{j_1, \beta} = \psi^{\ell'}_{\beta}$.  
Then $H^{(i_1)}_{j_1} + \Gamma_{i_1} = \ell' E$ for some $\Z$-divisor $E$, 
which contradicts (i) and (ii). 
Hence, $L_1$ is a field. 
We may apply the same argument for $L_2:=L_1[t_{i_2, j_2}]/(t^{\ell}_{i_2, j_2}-\varphi^{(i_2)}_{j_2, \alpha})$ 
and the normalisation $Y_1$ of $X$ in $L_1$, because 
$Y_1 \to X$ is \'etale outside $H^{(i_1)}_{j_1} + \Gamma_{i_1}$ and 
the pullback $(H^{(i_2)}_{j_2})_1$ of $H^{(i_2)}_{j_2}$ on $Y_1$ is reduced 
around the generic points of $(H^{(i_2)}_{j_2})_1$. 
Repeating this procedure finitely many times, we see that 
$K(X)[t_{i, j}]_{i, j}/(t_{i, j}^{\ell} - \varphi^{(i)}_{j, \alpha})_{i, j}$ is a field. 
Thus, (iii) holds. 

Fix $\alpha \in A$ and a closed point $x \in U_{\alpha} \subset X$. 
Set $I_x :=\{i \in I\,|\, x \in \Gamma_i\}$. 
For each $i \in I_x$, we fix $j_i$ such that $1 \leq j_i \leq d$ and $x \not\in H^{(i)}_{j_i}$, 
where the  existence of such $j_i$ is guaranteed 
by the choice of $\{H^{(i)}_j\}_{i, j}$. 
Set 
\begin{equation}\label{e1-Kawamata-cover}
\begin{split}
R_x := \{\varphi^{(i)}_{j_i \alpha}\,|\, i \in I_x\} 
&\cup \{\varphi^{(i)}_{j \alpha}\,|\, i \not\in I_x, x\in H^{(i)}_j\} \\
&\cup \{\varphi^{(i)}_{j\alpha}/\varphi^{(i)}_{j_i\alpha}\,|\, i \in I_x, x \in H^{(i)}_j\}
\end{split}
\end{equation}
\begin{equation}\label{e2-Kawamata-cover}
T_x:=\{\varphi^{(i)}_{j\alpha}/\varphi^{(i)}_{j_i\alpha}\,|\, i \in I_x, x \not\in H^{(i)}_j\}
\cup \{\varphi^{(i)}_{j\alpha}\,|\, i \not\in I_x, x \not\in H^{(i)}_{j}\}. 
\end{equation}
Note that $R_x$ forms a part of a regular system of parameters of $\MO_{X, x}$ and 
any element of $T_x$ is a unit of $\MO_{X, x}$. 
Indeed, the three sets in the right hand side of (\ref{e1-Kawamata-cover}) corresponds to the following sum: 
\[
\left( \sum_{i \in I_x} \Gamma_i \right) 
+ \left( \sum_{i \not\in I_x, x \in H^{(i)}_j} H^{(i)}_j\right) 
+ \left( \sum_{i \in I_x, x \in H^{(i)}_j} H^{(i)}_j\right). 
\]
We see that each element of $T_x$ is a unit of $\MO_{X, x}$ in a similar way. 
Hence, there is a regular system of parameters $z_1, ..., z_d$ of $\MO_{X, x}$ such that 
$R_x= \{z_1, ..., z_e\}$ with $0 \leq e \leq d$. 
For $T_x=\{u_1, ..., u_s\}$, 
we have that 
\[
K(Y) 
= K(X) [(\varphi^{(i)}_{j, \alpha})^{1/\ell}]_{i, j} =K(X)[z_1^{1/\ell}, ..., z_e^{1/\ell}, u_1^{1/\ell}, ..., u_s^{1/\ell}],
\]
where the second equality holds by (\ref{e1-Kawamata-cover}) and 
(\ref{e2-Kawamata-cover}). 
For $A:=\MO_{X, x}$, we set 
\[
B:=A[z_1^{1/\ell}, ..., z_e^{1/\ell}, u_1^{1/\ell}, ..., u_s^{1/\ell}]. 
\]
It is clear that $B$ is a finitely generated $A$-module. 
Furthermore, the function field of $\Spec\,B$ coincides with $K(Y)$. 
By (iii), the set 
\[
\{z_1^{c_1/\ell} \cdots z_e^{c_e/\ell} u_1^{d_1/\ell} \cdots u_s^{d_s/\ell}\,|\, 0 \leq c_i, d_j \leq \ell -1\}
\]
is linearly independent over $K(X)$. 
In particular, it holds that 
\[
B \simeq A[x_i, y_j]_{i, j}/(x_i^{\ell}-z_i, y_j^{\ell}-u_j)_{i, j}.
\]
It follows from Lemma \ref{l-Kawamata-cover} that $B$ is a regular ring. 
In particular, $B$ coincides with the integral closure of $A$ in $K(Y)$. 
Hence, $Y$ is a smooth variety and (1) holds.

It suffices to show (3). 
Fix a prime divisor $\Gamma_i$, a closed point $x \in \Gamma_i$, and 
$\alpha \in A$ with $x \in U_{\alpha}$. 
Set $V_{\alpha} := U_{\alpha} \setminus H_{j_i}^{(i)}$, which is an open neighbourhood of $x \in U_{\alpha}$. 
By $
\Gamma_i|_{V_{\alpha}}=
(H^{(i)}_{j_i} + \Gamma_i)|_{V_{\alpha}} = {\rm div}(\varphi^{(i)}_{j_i, \alpha})|_{V_{\alpha}}$, 
we obtain
\[
g^*(\Gamma_i)|_{g^{-1}(V_{\alpha})} 
= 
g^*({\rm div}(\varphi^{(i)}_{j_i, \alpha}))|_{g^{-1}(V_{\alpha})} = 
\ell {\rm div}( (\varphi^{(i)}_{j_i, \alpha})^{1/\ell})|_{g^{-1}(V_{\alpha})}, 
\]
which implies that $g^*(\Gamma_i) = \ell \Gamma'_i$ for some effective $\Z$-divisor $\Gamma'_i$ on $Y$. Thus, (3) holds. 
\end{proof}

\begin{prop}\label{p-Kawamata}
Let $k$ be a perfect field of characteristic $p>0$.  
Let $X$ be a quasi-projective smooth geometrically connected variety over $k$. 
Fix $\ell \in \Z_{>0} \setminus p\Z$. 
Let $D$ be a $\Q$-divisor such that 
$\{D\}$ is simple normal crossing and $\ell D$ is a $\Z$-divisor. 
Then there exist finite surjective $k$-morphisms 
\[
h:X'' \xrightarrow{g} X' \xrightarrow{f} X
\]
of smooth $k$-varieties that satisfy the following properties. 
\begin{enumerate}
\item 
There exists a finite extension $k \subset k'$ of perfect fields 
such that $X'= X \times_k k'$, $f:X' \to X$ is the projection, 
and $f$ is weakly $\ell'$-cyclic for some $\ell' \in \Z_{>0} \setminus p\Z$. 
\item 
$g$ is weakly $\ell$-cyclic. 
\item  $h^*D$ is a $\Z$-divisor. 
\end{enumerate}
\end{prop}

\begin{proof}
Fix an algebraic closure $\overline k$ of $k$. 
Let $\overline{\mathbb F}_p$ be the algebraic closure of $\mathbb F_p$ 
that is a subfield of $\overline k$. 
Let $k_1$ be the smallest subfield of $\overline k$ that contains 
$\overline{\mathbb F}_p \cup k$.  

Set $X_1:=X \times_k k_1$ and we apply Theorem \ref{t-Kawamata-cover}. 
Then there exists a finite surjective $k_1$-morphism 
\[
g_1:Y_1 \to X_1
\]
which satisfies the properties (1)--(3) listed in Theorem \ref{t-Kawamata-cover}. 
We can find an intermediate field $k \subset k' \subset k_1$ and 
a finite surjective morphism $g:X'' \to X'$ of smooth $k'$-varieties 
such that $[k':k] < \infty$ and $g \times_{k'} k_1 = g_1$. 

It is enough to show that $f$ is weakly $\ell'$-cyclic for some $\ell' \in \Z_{>0} \setminus p\Z$. 
Since $k_1$ is the composite field of $\F_{p^{e}}$ and $k$ for some $e \in \Z_{>0}$, 
it holds that $k_1=k(\zeta_{\ell'})$ for a primitive $\ell'$-th root $\zeta_{\ell'}$ of the unity, 
where $\ell':=p^{e}-1$. 
Then $f$ is weakly $\ell'$-cyclic. 
\end{proof}

\begin{thm}\label{t-KVV}
Let $k$ be a perfect field of characteristic $p>0$. 
Let $f:X \to Y$ be a projective $k$-morphism 
from an $N$-dimensional smooth $k$-variety $X$ to 
a scheme $Y$ of finite type over $k$. 
Let $A$ be an $f$-ample $\Z_{(p)}$-divisor on $X$, 
i.e. there exists $\ell \in \Z_{>0} \setminus p\Z$ such that 
$\ell A$ is an $f$-ample $\Z$-divisor on $X$. 
Then 
\[
R^if_*(\mathcal Hom_{W\MO_X}(W\MO_X(-A), W\Omega_X^N)=0
\]
for any $i \in \Z_{>0}$. 
\end{thm}

\begin{proof}
Replacing $f:X \to Y$ by $f':X \to Y'$ for the Stein factorisation $f:X \xrightarrow{f'} Y' \to Y$ of $f$, 
we may assume that $f_*\MO_X=\MO_Y$. 
Furthermore, taking the algebraic closure of $k$ in $K(X)$, 
the problem is reduced to the case when $X$ is geometrically connected over $k$. 
Fix $\ell \in \Z_{>0} \setminus p\Z$ such that $\ell A$ is a $\Z$-divisor. 
Let 
\[
h:X'' \xrightarrow{g} X' \xrightarrow{f} X
\]
be as in Proposition \ref{p-Kawamata}. 
For $V \in \{X, X', X''\}$, 
let $f_V: V \to Y$ be the induced morphism and 
set $A_V$ to be the pullback of $A$ to $V$. 
Let $\varphi:Z' \to Z$ be one of $g$ and $f$.

\setcounter{step}{0}
\begin{step}\label{s1-KVV}
The $W\MO_Z$-module homomorphism 
\[
\mathcal Hom_{W\MO_Z}(\varphi_*W\MO_{Z'}(-A_{Z'}), W\Omega_Z^N) 
\to 
\mathcal Hom_{W\MO_Z}(W\MO_Z(-A_Z), W\Omega_Z^N),
\]
induced by $W\MO_Z(-A_Z) \to \varphi_*W\MO_{Z'}(-A_{Z'})$, is a split surjection. 
\end{step}

\begin{proof}[Proof of Step \ref{s1-KVV}]
We have the splitting: 
\[
{\rm id}: W\MO_Z(-A_Z) \xrightarrow{\alpha} \varphi_*(W\MO_{Z'}(-A_{Z'})) \xrightarrow{\beta} W\MO_Z(-A_Z). 
\]
Indeed, if $\varphi=f$ (resp. $\varphi=g$), 
then the above splitting holds 
by Proposition \ref{p-pull-field-ext}(3) (resp. Lemma \ref{l-finite-trace}(2)). 
Applying the contravariant functor $\mathcal Hom_{W\MO_Z}(-, W\Omega_Z^N)$, 
we obtain the required splitting surjection. 
This completes the proof of Step \ref{s1-KVV}. 
\end{proof}

\begin{step}\label{s2-KVV}
For any $n \in \Z_{>0}$, 
the induced $W\MO_Z$-module homomorphism 
\[
\mathcal Hom_{W_n\MO_Z}(\varphi_*W_n\MO_{Z'}(-A_{Z'}), W_n\Omega_Z^N)
\to 
\mathcal Hom_{W\MO_Z}(\varphi_*W\MO_{Z'}(-A_{Z'}), W_n\Omega_Z^N)
\]
is an isomorphism. 
\end{step}

\begin{proof}[Proof of Step \ref{s2-KVV}]
Since $\varphi:Z' \to Z$ is weakly $\ell''$-cyclic for some $\ell'' \in \Z_{>0} \setminus p\Z$, the assertion of Step \ref{s2-KVV} follows from Proposition \ref{p-WO-vs-WO-Y} 
and Lemma \ref{l-DRW-inje}. 
\end{proof}

\begin{step}\label{s3-KVV}
There exists an isomorphism of  $W\MO_Z$-modules: 
\[
\mathcal Hom_{W\MO_Z}(\varphi_*W\MO_{Z'}(-A_{Z'}), W\Omega_Z^N)) 
\simeq 
\varphi_* \mathcal Hom_{W\MO_{Z'}}(W\MO_{Z'}(-A_{Z'}), W\Omega_{Z'}^N). 
\]
\end{step}

\begin{proof}[Proof of Step \ref{s3-KVV}]
It holds that  
\[
\mathcal Hom_{W\MO_Z}(\varphi_*W\MO_{Z'}(-A_{Z'}), W\Omega_Z^N)\\
\simeq 
\varprojlim_n (\mathcal Hom_{W\MO_Z}(\varphi_*W\MO_{Z'}(-A_{Z'}), W_n\Omega_Z^N). 
\]
We have the following isomorphisms, each of which is canonically defined: 
\begin{eqnarray*}
&& \mathcal Hom_{W\MO_Z}(\varphi_*W\MO_{Z'}(-A_{Z'}), W_n\Omega_Z^N) \\
&\simeq & 
\mathcal Hom_{W_n\MO_Z}(\varphi_*W_n\MO_{Z'}(-A_{Z'}), W_n\Omega_Z^N)\\
&\simeq &
\varphi_*\mathcal Hom_{W_n\MO_{Z'}}(W_n\MO_{Z'}(-A_{Z'}), \varphi^!W_n\Omega_Z^N)\\
&\simeq &
\varphi_*\mathcal Hom_{W_n\MO_{Z'}}(W_n\MO_{Z'}(-A_{Z'}), W_n\Omega_{Z'}^N)\\
&\simeq &
\varphi_*\mathcal Hom_{W\MO_{Z'}}(W\MO_{Z'}(-A_{Z'}), W_n\Omega_{Z'}^N),\\
\end{eqnarray*}
where the first isomorphism holds  by Step \ref{s2-KVV}, 
the second one follows from the Grothendieck duality, 
and the third one holds by \cite[Theorem 4.1]{Eke84}. 
This completes the proof of Step \ref{s3-KVV}. 
\end{proof}

\begin{step}\label{s4-KVV}
For $V \in \{X, X', X''\}$, 
consider the following assertion $(V)$. 
\begin{enumerate}
\item[$(V)$] 
Let $g:V \to Y'$ be a projective $Y$-morphism to a projective $Y$-scheme $Y'$. 
Then the equation 
\[
R^ig_*(\mathcal Hom_{W\MO_{V}}(W\MO_{V}(-A_{V}), W\Omega_{V}^N)=0
\]
holds for any $i \in \Z_{>0}$. 
\end{enumerate}
Then the following holds. 
\begin{enumerate}
\item If $(Z')$ holds, then $(Z)$ holds. 
\item $(X'')$ holds. 
\item $(X)$ holds. 
\end{enumerate}
\end{step}

\begin{proof}[Proof of Step \ref{s4-KVV}]
Let us show (1). 
Assume $(Z')$. 
Let $g: Z \to Y'$ be as in the statement. 
We have the following isomorphisms in the derived category of $W\MO_{Y'}$-modules: 
\begin{eqnarray*}
&&Rg_*\varphi_* \mathcal Hom_{W\MO_{Z'}}(W\MO_{Z'}(-A_{Z'}), W\Omega_{Z'}^N)\\
&\simeq &
Rg_*R\varphi_* \mathcal Hom_{W\MO_{Z'}}(W\MO_{Z'}(-A_{Z'}), W\Omega_{Z'}^N)\\
&\simeq &
R(g \circ \varphi)_* \mathcal Hom_{W\MO_{Z'}}(W\MO_{Z'}(-A_{Z'}), W\Omega_{Z'}^N)\\
&\simeq &
(g \circ \varphi)_* \mathcal Hom_{W\MO_{Z'}}(W\MO_{Z'}(-A_{Z'}), W\Omega_{Z'}^N)\\
\end{eqnarray*}
where the first and third isomorphisms hold by $(Z')$. 
In particular, we get 
\[
R^ig_*\varphi_* \mathcal Hom_{W\MO_{Z'}}(W\MO_{Z'}(-A_{Z'}), W\Omega_{Z'}^N) =0
\]
for any $i>0$.
 By Step \ref{s3-KVV}, the equation 
\[
R^ig_* \mathcal Hom_{W\MO_{Z}}(\varphi_*W\MO_{Z'}(-A_{Z'}), W\Omega_{Z}^N)) = 0
\]
holds for any $i>0$. 
Then Step \ref{s1-KVV} implies $(Z)$. 
Thus, (1) holds. 

The assertion (2) follows from Theorem \ref{t-Kod-van}(4) and 
the fact that $A_{X''}$ is a $\Z$-divisor. 
The assertion (3) holds by (1) and (2). 
This completes the proof of Step \ref{s4-KVV}. 
\end{proof}
Step \ref{s4-KVV}(3) completes the proof of Theorem \ref{t-KVV}. 
\end{proof}

\begin{cor}\label{c-KVV}
Let $k$ be a perfect field of characteristic $p>0$. 
Let $f:X \to Y$ be a projective $k$-morphism 
from an $N$-dimensional smooth $k$-variety $X$ to 
a scheme $Y$ of finite type over $k$. 
Let $A$ be an $f$-ample $\Q$-divisor on $X$. 
Then
\[
R^if_*(\mathcal Hom_{W\MO_X}(W\MO_X(-A), W\Omega_X^N)_{\Q}=0
\]
for any $i \in \Z_{>0}$. 
\end{cor}

\begin{proof}
In the following argument, we denote by 
\begin{equation}\label{e1-(4)-KVV}
F_X:(X, W\MO_{X, \Q}) \xrightarrow{\simeq} (X, W\MO_{X, \Q}), 
\end{equation}
\[
F_Y:(Y, W\MO_{Y, \Q}) \xrightarrow{\simeq} (Y, W\MO_{Y, \Q})
\]
the isomorphisms of ringed spaces induced by Frobenius. 
Note that these morphisms are isomorphisms by $p=FV = VF$. 
In particular, we have that 
\[
(F_X)_* \circ (F_X)^*=(F_X)^* \circ (F_X)_*={\rm id}.
\]
Fix a positive integer $e$ such that $p^eA$ is a $\Z_{(p)}$-divisor. 
Then we obtain isomorphisms of the derived category of $W\MO_{X, \Q}$-modules 
(cf. Subsection \ref{ss-notation}(\ref{ss-notation-Q})): 
\begin{eqnarray*}
&& Rf_*\mathcal Hom_{W\MO_{X, \Q}}(W\MO_X(-A)_{\Q}, W\Omega^N_{X, \Q})\\
&\simeq & Rf_*(F_X^e)_*(F_X^e)^*\mathcal Hom_{W\MO_{X, \Q}}(W\MO_X(-A)_{\Q}, W\Omega^N_{X, \Q}))\\
&\simeq &  (F_Y^e)_*Rf_*\mathcal Hom_{(F_X^e)^*W\MO_{X, \Q}}((F_X^e)^*W\MO_X(-A)_{\Q}, (F_X^e)^*W\Omega^N_{X, \Q}))\\
&\simeq &  (F_Y^e)_*Rf_*\mathcal Hom_{W\MO_{X, \Q}}(W\MO_X(-p^e A)_{\Q}, 
(F_X^e)^*(F_X^e)_*W\Omega^N_{X, \Q}))\\
&\simeq &  (F_Y^e)_*Rf_*\mathcal Hom_{W\MO_{X, \Q}}(W\MO_X(-p^e A)_{\Q}, W\Omega^N_{X, \Q}))\\
&\simeq &(F_Y^e)_*f_*\mathcal Hom_{W\MO_{X, \Q}}(W\MO_X(-p^e A)_{\Q}, W\Omega^N_{X, \Q})),
\end{eqnarray*}
where the first and fourth isomorphisms hold by (\ref{e1-(4)-KVV}), 
the second one by the fact that $(F_X^e)_*$ 
and $(F_Y^e)_*$ are exact functors, 
the third one by Theorem \ref{t-descent}, and 
the last one by Theorem \ref{t-KVV}
\end{proof}

\begin{thm}\label{t-KVV-minus}
Let $k$ be a perfect field of characteristic $p>0$ and 
let $X$ be an $N$-dimensional projective Cohen-Macaulay normal variety over $k$. 
Let $A$ be an ample $\R$-Cartier $\R$-divisor on $X$. 
Then the following hold. 
\begin{enumerate}
\item 
There exists a positive integer $s_1$ such that 
$$H^j(X, W_n\MO_X(-sA))=0$$ 
for any $n  \in \Z_{>0}$, $j \in \Z_{<N}$ and $s \in \R_{\geq s_1}$. 
\item 
There exists a positive integer $s_2$ such that 
$$H^j(X, W\MO_X(-sA))=0$$ 
for any $j \in \Z_{<N}$ and $s \in \R_{\geq s_2}$. 
\item There exists a positive integer $t$ such that 
\[
H^j(X, W\MO_X(-A))
\]
is $p^t$-torsion for any $j \in \Z_{<N}$. 
In particular, the equation 
\[
H^j(X, W\MO_X(-A))_{\Q}=0 
\]
holds for any $j \in \Z_{<N}$. 
\end{enumerate}
\end{thm}

\begin{proof}
For any $\R$-divisor $D$ on $X$, 
Proposition \ref{p-div-induction} induces an exact sequence of $W\MO_X$-homomorphisms: 
$$0 \to (F_X)_*(W_n\MO_X(-psD)) \xrightarrow{V} W_{n+1}\MO_X(-sD) \to \MO_X(-sD) \to 0.$$
Thus, by the Serre vanishing theorem \cite[Ch. III, Theorem 7.6]{Har77}, there exists a positive integer $s_1$ such that 
$$H^j(X, W_n\MO_X(-sA))=0$$
for any $n  \in \Z_{>0}$, $j \in \Z_{<N}$ and $s \in \R_{\geq s_1}$. 
Thus (1) holds. Set $s_2:=s_1$. 
Since $H^j(X, W_n\MO_X(-sA))$ is a finitely generated $W_n(k)$-module, 
$H^j(X, W_n\MO_X(-sA))$ is an artinian $W(k)$-module. 
Therefore, the projective system $\{H^j(X, W_n\MO_X(-sA))\}_{n \in \Z_{>0}}$ satisfies the Mittag--Leffler condition. It follows from Lemma \ref{l-CR12}(3) that the isomorphism 
$$H^j(X, W\MO_X(-sA)) \simeq \varprojlim_n H^j(X, W_n\MO_X(-sA))=0$$ 
holds 
for any $j<N$ and $s \in \R_{\geq s_2}$. 
Hence (2) holds. 

Let us show (3). 
Fix $t \in \Z_{>0}$ such that $p^t \geq s_2$. 
The multiplication map $p^t$ have the following factorisation: 
\[
p^t: W\MO_X(-A) \xrightarrow{F^t} (F^t_X)_*(W\MO_X(-p^tA)) \xrightarrow{V^t} W\MO_X(-A). 
\]
Take cohomologlies $H^j(X, -)$: 
{\small 
\[
p^t:
H^j(X, W\MO_X(-A)) \xrightarrow{F^t} H^j(X, (F^t_X)_*(W\MO_X(-p^tA))) \xrightarrow{V^t} 
H^j(X, W\MO_X(-A)).
\]
}
\noindent
The middle term $H^j(X, (F^t_X)_*(W\MO_X(-p^tA)))$ is zero by the following computation: 
\[
H^j(X, (F^t_X)_*(W\MO_X(-p^tA))) \simeq H^j(X, W\MO_X(-p^tA)) = 0, 
\]
where the isomorphism holds by the fact that $(F^t_X)_*$ is exact and 
the equality follows from the choice of $t$. 
Therefore, $H^j(X, W\MO_X(-A))$ is $p^t$-torstion, as desired. 
\end{proof}

\begin{rem}\label{r-KVV-minus}
Let $k$ be a perfect field of characteristic $p>0$ and 
let $X$ be a proper scheme over $k$. 
Assume that $H$ is a closed subscheme such that $U = X \setminus H$ 
is affine, smooth, and pure $N$-dimensional. 
Then, Berthelot--Bloch--Esnault proved that the equation 
\[
H^j(X, WI_H)_{\Q}=0
\]
holds for any $j \neq N$ (\cite[Corollary 1.2]{BBE07}). 
If $X$ is smooth and $H$ $H$ is an ample effective Cartier divisor, 
this result is very similar to but different from Theorem \ref{t-KVV-minus}. 
This is because our Witt divisorial sheaf $W\MO_X(-H)$ differs from the Witt ideal sheaf $WI_H$. 
\end{rem}


\section{Examples}\label{s-examples}

\subsection{Ample divisors}

In this subsection, we compute some cohomologies 
for Teichm\"uller lifts of ample invertible sheaves. 
First we give a Witt analogue of the Serre vanishing theorem (Theorem \ref{t-Witt-Serre}). 
Second we compute dimension of cohomologies for Teichm\"uller lifts of ample invertible sheaves (Theorem \ref{t-ample-dim}). 

\begin{thm}\label{t-Witt-Serre}
Let $f:X \to Y$ be a projective morphism of $F$-finite noetherian $\F_p$-schemes. 
Let $A$ be an $f$-ample invertible sheaf on $X$ and 
let $I$ be a coherent ideal sheaf on $X$. 
Then the equation 
$$R^if_*(WI \otimes_{W\MO_X} \underline{A})_{\Q}=0$$
holds for any positive integer $i>0$. 
\end{thm}

\begin{proof}
We can find a positive integer $s_1$ such that 
\begin{equation}\label{e-Witt-Serre1}
R^if_*(I \otimes_{\MO_X} A^{\otimes s})=0
\end{equation}
for any $i \in \Z_{>0}$ and any $s \in \Z_{\geq s_1}$. 
We have an exact sequence of $W\MO_X$-modules: 
\begin{equation}\label{e-Witt-Serre2}
0 \to (F_X)_*(W_nI) \xrightarrow{V} W_{n+1}I \to I \to 0.
\end{equation}
By (\ref{e-Witt-Serre1}) and (\ref{e-Witt-Serre2}), we get  
\begin{equation}\label{e-Witt-Serre3}
R^if_*(W_nI \otimes_{W\MO_X} {\underline{A}}^{\otimes s})=0
\end{equation}
for any $i \in \Z_{>0}$, $s \in \Z_{\geq s_1}$, and $n\in \Z_{>0}$. 
Fix $t \in \Z_{>0}$ such that $p^t \geq s_1$. 
For 
\[
M:=WI \otimes_{W\MO_X} {\underline{A}}^{\otimes p^t} \quad 
{\rm and} \quad 
M_n:=W_nI \otimes_{W\MO_X} {\underline{A}}^{\otimes p^t},
\]
the induced homomorphism $R : M_{n+1} \to M_n$ is surjective, hence 
it follows from  Lemma \ref{l-CR12}(2) that 
$R\varprojlim_n M_n \simeq \varprojlim_n M_n \simeq  M$. 
Then we obtain isomorphisms 
{\small
$$Rf_*M \simeq \left(Rf_* \circ R\varprojlim_n\right) M_n 
\simeq \left(R\varprojlim_n \circ Rf_*\right) M_n
\simeq R\varprojlim_n (f_*M_n) 
\simeq \varprojlim_n (f_*M_n),$$
}
where the third isomorphism follows from (\ref{e-Witt-Serre3}) 
and the fourth isomorphism holds by Lemma \ref{l-CR12}(2)  
and the fact that an exact sequence 
$$0 \to (F_X^n)_*I \to W_{n+1}I \xrightarrow{R} W_nI \to 0$$
implies that the projective system $\{f_*M_n\}_{n \in \Z_{>0}}$ consists of 
surjective homomorphisms $f_*(R \otimes_{W\MO_X} \underline{A}^{\otimes p^t})$. 
Therefore, we have that 
\begin{equation}\label{e-Witt-Serre4}
R^if_*(WI \otimes_{W\MO_X} {\underline{A}}^{\otimes p^t}) = R^if_*M=0
\end{equation}
for any $i>0$. 
By $p^t = F^tV^t =V^tF^t$, the $W\MO_{X, \Q}$-module homomorphism 
\[
F^t: WI_{\Q} \xrightarrow{\simeq}  (F^t_X)_*(WI_{\Q})
\]
is an isomorphism. 
Then it holds that 
\begin{eqnarray*}
R^if_*(WI \otimes_{W\MO_X} \underline{A})_{\Q}
&\simeq& R^if_*((F^t_X)_*(WI) \otimes_{W\MO_X} \underline{A})_{\Q}\\
&\simeq& R^if_*((F^t_X)_*(WI \otimes_{W\MO_X} \underline{A}^{p^t}))_{\Q}\\
&\simeq &(F^t_Y)_*R^if_*(WI \otimes_{W\MO_X} \underline{A}^{p^t})_{\Q}\\
&=&0,
\end{eqnarray*}
where the second isomorphism follows from the projection formula, 
the third one holds by the fact that the functors $(F^t_X)_*$ and $(F^t_Y)_*$ 
are exact, and the last equation follows from (\ref{e-Witt-Serre4}). 
\end{proof}

\begin{thm}\label{t-ample-dim}
Let $k$ be a perfect field of characteristic $p>0$ and set $Q:=W(k)_{\Q}$. 
Let $X$ be a projective scheme over $k$ such that $\dim X>0$. 
Let $A$ be an ample invertible sheaf on $X$. 
Then the following hold. 
\begin{enumerate}
\item 
$H^i(X, \underline{A})_{\Q}=0$ for $i>0$. 
\item 
$\dim_Q H^0(X, \underline{A})_{\Q}=\infty$. 
\end{enumerate}
\end{thm}

\begin{proof}
Since (1) follows directly from Theorem \ref{t-Witt-Serre}, 
let us prove (2). 
Fix an arbitrary positive integer $m$. 
By $\dim X>0$, we can find distinct closed points $P_1, \cdots, P_m$ of $X$. 
Let $Y$ be the reduced zero-dimensional closed subscheme of $X$ 
that is set-theoretically equal to $\{P_1, \cdots, P_m\}$. 
By the exact sequence 
$$0 \to WI_Y \otimes_{W\MO_X} \underline{A} \to \underline{A} \to \underline{A}|_Y \to 0,$$
Theorem \ref{t-Witt-Serre} implies that the induced map 
$$H^0(X, \underline{A})_{\Q} \to H^0(Y, \underline{A}|_Y)_{\Q}$$
is surjective. Hence, it holds that  
$$\dim_Q H^0(X, \underline{A})_{\Q} \geq 
\dim_Q H^0(Y, \underline{A}|_Y)_{\Q}=
\dim_Q H^0(Y, W\MO_Y)_{\Q} \geq m.$$
Since $m$ was chosen to be an arbitrary positive integer, 
the assertion (2) holds. 
\end{proof}

\subsection{Anti-ample divisors}

\subsubsection{The highest cohomology}

The purpose of this subsection is to prove Theorem \ref{t-highest-infinite}, 
which can be considered as a dual of Theorem \ref{t-ample-dim}(2).

\begin{thm}\label{t-highest-infinite}
Let $k$ be a perfect field of characteristic $p>0$ and set $Q:=W(k)_{\Q}$. 
Let $X$ be an $N$-dimensional smooth projective scheme over $k$ 
such that $N>0$. 
Let $A$ be an ample invertible sheaf on $X$. 
Then $\dim_Q H^N(X, \underline{A}^{-1})_{\Q}=\infty$. 
\end{thm}

\begin{proof}
Replacing $X$ by an $N$-dimensional irreducible component, we may assume that $X$ is connected. 

\setcounter{step}{0}

\begin{step}\label{s1-highest-infinite}
There exists a positive integer $s_1$ such that the map 
$$F:H^N(X, A^{-s}) \to H^N(X, (F_X)_*(A^{-sp}))$$
induced by Frobenius $\MO_X \to (F_X)_*\MO_X$ is injective for any $s \in \Z_{\geq s_1}$. 
\end{step}

\begin{proof}[Proof of Step \ref{s1-highest-infinite}]
We have an exact sequence 
$$0 \to \MO_X \to (F_X)_*\MO_X \to B \to 0$$
for a coherent locally free $\MO_X$-module $B$ (Proposition \ref{p-BZ}). 
By \cite[Ch. III, Theorem 7.6(b)]{Har77}, we have that 
$$H^{N-1}(X, B \otimes_{\MO_X} A^{-s})=0$$
for $s \gg 0$. 
This completes the proof of Step \ref{s1-highest-infinite}.  
\end{proof}

\begin{step}\label{s2-highest-infinite}
There exists a positive integer $s_2 \in \Z_{>0}$ such that the map 
$$F:H^N(X, \underline{A}_{\leq n}^{-s}) \to H^N(X, (F_X)_*(\underline{A}_{\leq n}^{-sp}))$$
induced by Frobenius $W_n\MO_X \to (F_X)_*(W_n\MO_X)$ is injective 
for any $s \in \Z_{\geq s_2}$ and $n \in \Z_{>0}$. 
\end{step}

\begin{proof}[Proof of Step \ref{s2-highest-infinite}] 
For any invertible sheaf $L$ on $X$, we have a commutative diagram 
of $W\MO_X$-homomorphisms: 
$$\begin{CD}
0 @>>> (F_X)_*(\underline{L}^p_{\leq n}) @>V >> \underline{L}_{\leq n+1} @>>> L @>>> 0\\
@. @VV(F_X)_*F V @VVFV @VVFV\\
0 @>>> (F^2_X)_*(\underline{L}_{\leq n}^{p^2}) @>V >> (F_X)_*(\underline{L}_{\leq n+1}^p) @>>> (F_X)_*(L^{p}) @>>> 0\\
\end{CD}$$
where both the horizontal sequences are exact (Proposition \ref{p-div-induction}). 
For any $s \gg 0$, it holds that $H^{N-1}(X, A^{-s})=0$, 
hence we get a commutative diagram: 
{\tiny 
$$\begin{CD}
0 @>>> H^N((F_X)_*(\underline{A}^{-sp}_{\leq n})) @>>> H^N(\underline{A}^{-s}_{\leq n+1}) @>>> H^N(A^{-s}) @>>> 0\\
@. @VV(F_X)_*F V @VVF V @VVF V\\
0 @>>> H^N((F^2_X)_*(\underline{A}^{-sp^2}_{\leq n})) @>>> H^N((F_X)_*(\underline{A}^{-sp}_{\leq n+1})) @>>> H^N((F_X)_*(A^{-sp})) @>>> 0,\\
\end{CD}$$
}
where both the horizontal sequences are exact. 
By the snake lemma and induction on $n$, 
Step \ref{s1-highest-infinite} implies Step \ref{s2-highest-infinite}. 
\end{proof}

\begin{step}\label{s3-highest-infinite}
There exists a positive integer $s_3 \in \Z_{>0}$ such that the map 
$$F:H^N(X, \underline{A}^{-s}) \to H^N(X, (F_X)_*(\underline{A}^{-sp}))$$
induced by Frobenius $W\MO_X \to (F_X)_*(W\MO_X)$ is injective for any 
$s \in \Z_{\geq s_3}$. 
\end{step}

\begin{proof}[Proof of Step \ref{s3-highest-infinite}]
Set $s_3:=s_2$. 
Fix $i \in \Z_{\geq 0}, e \in \Z_{\geq 0},$ and an invertible $\MO_X$-module $L$. 
Since $H^i(X, (F_X^e)_*\underline{L}_{\leq n})$ is a finitely generated $W_n(k)$-module for any $i \in \Z_{\geq 0}$, 
$H^i(X, (F_X^e)_*\underline{L}_{\leq n})$ is an artinian $W(k)$-module. 
Therefore, the projective system 
\[
\{H^i(X, (F_X^e)_*\underline{L}_{\leq n})\}_{n \in \Z_{>0}}
\]
satisfies the Mittag--Leffler condition. 
Hence,  
Lemma \ref{l-CR12}(3) implies that 
$$H^i(X, (F_X^e)_*\underline{L}) \simeq \varprojlim_n H^i(X, (F_X^e)_*\underline{L}_{\leq n}).$$
Thus Step \ref{s3-highest-infinite} follows from Step \ref{s2-highest-infinite}. 
\end{proof}

\begin{step}\label{s4-highest-infinite}
There exists a positive integer $s_4 \in \Z_{>0}$ such that the map 
$$V:H^N(X, (F_X)_*(\underline{A}^{-sp})) \to H^N(X, \underline{A}^{-s})$$
induced by Verschiebung $V:(F_X)_*(W\MO_X) \to W\MO_X$ is injective for any $s \in \Z_{\geq s_4}$. 
\end{step}

\begin{proof}[Proof of Step \ref{s4-highest-infinite}]
The assertion follows from 
the exact sequence (Proposition \ref{p-div-induction}):  
$$0 \to (F_X)_*(\underline{A}^{-sp}) \xrightarrow{V} \underline{A}^{-s} \to A^{-s} \to 0$$
and $H^{N-1}(X, A^{-s})=0$ for any $s \gg 0$. 
\end{proof}

\begin{step}\label{s5-highest-infinite}
There exists a positive integer $s_5 \in \Z_{>0}$ such that 
$H^N(X, \underline{A}^{-s})$ is a $p$-torsion free $W(k)$-module 
for any $s \in \Z_{\geq s_5}$. 
\end{step}

\begin{proof}[Proof of Step \ref{s5-highest-infinite}] 
Set $s_5:=\max\{s_3, s_4\}$. 
Since $p=FV=VF$, the assertion follows from Step \ref{s3-highest-infinite} 
and Step \ref{s4-highest-infinite}. 
\end{proof}

\begin{step}\label{s6-highest-infinite}
Let $\nu_0$ be a positive integer. 
Then there exists a positive integer $s_0 \in \Z_{>0}$ such that 
for any $s \in \Z_{\geq s_0}$, $H^N(X, \underline{A}^{-s})$ 
is a $p$-torsion free $W(k)$-module whose rank is larger than $\nu_0$. 
\end{step}

\begin{proof}[Proof of Step \ref{s6-highest-infinite}]
By Serre duality, we can find a positive integer $s'_0$ such that 
$$\dim_k H^N(X, A^{-s}) \geq \nu_0+1$$ 
for any $s \in \Z_{\geq s'_0}$. 
Set $s_0:=\max\{s'_0, s_5\}$. 
Fix $s \in \Z_{\geq s_0}$. 
By Step \ref{s5-highest-infinite}, 
we have that $H^N(X, \underline{A}^{-s})$ is a $p$-torsion free $W(k)$-module. 
Since the map 
$$R':H^N(X, \underline{A}^{-s}) \to H^N(X, A^{-s})$$
induced by $R:W\MO_X \to \MO_X$ is surjective, 
there exist elements $\zeta_0, \cdots, \zeta_{\nu_0} \in H^N(X, \underline{A}^{-s})$ 
such that 
$$R'(\zeta_0), \cdots, R'(\zeta_{\nu_0})$$
are linearly independent over $k$. 

It suffices to prove that $\zeta_0, \cdots, \zeta_{\nu_0}$ 
is linearly independent over $W(k)$. 
Assume that $\sum_{i=0}^{\nu_0}a_i\zeta_i=0$ 
for some $a_0, \cdots, a_{\nu_0} \in W(k)$. 
Then we have that $\sum_{i=0}^{\nu_0}\overline{a}_i R'(\zeta_i)=0$, 
where $\overline{a}_i :=a_i \mod pW(k)$. 
Since $R'(\zeta_0), \cdots, R'(\zeta_{\nu_0})$ are linearly independent over $k$, 
it holds that $\overline{a}_0=\cdots=\overline{a}_{\nu_0}=0$. 
Thus, we can write $a_i=pb_i$ for some $b_i \in W(k)$. 
Since $H^N(X, \underline{A}^{-s})$ is $p$-torsion free, 
we have that $\sum_{i=0}^{\nu_0}b_i\zeta_i=0$. 
Applying this argument again, 
we can find $c_0, \cdots, c_{\nu_0} \in W(k)$ such that $b_i=pc_i$ for any $i \in \{0, \cdots, \nu_0\}$. 
Hence, we see that $a_i \in p^2W(k)$ for any $i \in \{0, \cdots, \nu_0\}$. 
Repeating this procedure, we have that 
$$\{a_0. \cdots, a_{\nu_0}\} \subset \bigcap_{e=1}^{\infty} p^eW(k)=0.$$
Therefore, $\zeta_0, \cdots, \zeta_{\nu_0}$ 
is linearly independent over $W(k)$. 
This completes the proof of Step \ref{s6-highest-infinite}. 
\end{proof}

The assertion of Theorem \ref{t-highest-infinite} 
follows from Step \ref{s6-highest-infinite} and 
the isomorphisms of $Q$-vector spaces: 
\[
H^N(X, \underline{A}^{-1})_{\Q} \xrightarrow{F^t, \,\, \simeq} 
H^N(X, (F_X^t)_*\underline{A}^{-p^t})_{\Q} \simeq 
H^N(X, \underline{A}^{-p^t})_{\Q},
\]
where the first isomorphism holds because 
the $W\MO_{X, \Q}$-module homomorphism 
$F^t : \underline{A}^{-1}_{\Q} \to (F_X^t)_*\underline{A}^{-p^t}_{\Q}$ is an isomorphism by $p^t = F^t V^t =V^tF^t$. 
\end{proof}

\subsubsection{Cohomologies with torsion elements}

The following proposition shows that 
we actually need to take the tensor product with $\Q$ in the statement (2) of Theorem \ref{intro-main1}. 

\begin{prop}\label{p-cex-torsion}
Let $k$ be an algebraically closed field of characteristic $p>0$ 
and let $N$ be an integer satisfying $N \geq 2$. 
Then there exist 
an $N$-dimensional smooth projective variety $X$ over $k$ and
an ample invertible sheaf $A$ on $X$ such that 
$H^1(X, \underline{A}^{-1})\neq 0$ or $H^2(X, \underline{A}^{-1})\neq 0$.
\end{prop}

\begin{proof}
By \cite[Theorem 2]{Muk13}, 
there exist an $N$-dimensional smooth projective variety $X$ and 
an ample invertible sheaf $L$ on $X$ such that $H^1(X, L^{-1}) \neq 0$. 
By an exact sequence (Proposition \ref{p-div-induction}): 
$$0 \to (F_X)_*(\underline{L}^{-p}) \xrightarrow{V} \underline{L}^{-1} \to L^{-1} \to 0,$$
we obtain 
$$H^1(X, \underline{L}^{-1}) \neq 0 \quad {\rm or} \quad 
H^2(X, \underline{L}^{-p}) \neq 0.$$
Setting $A:=L$ or $A:=L^p$, we are done. 
\end{proof}


\begin{thebibliography}{00000000}



\bibitem[BBE07]{BBE07}
{P. Berthelot, S. Bloch, H. Esnault},
{\em On Witt vector cohomology for singular varieties,}
{Compos. Math. {\textbf{143}} (2007), no. 2, 363--392. }


\bibitem[BH93]{BH93}
{W. Bruns, J. Herzogs},
{\em Cohen-Macaulay rings,}
{Cambridge Studies in Advanced Mathematics, {\textbf{39}}. 
Cambridge University Press, Cambridge, 1993.}


\bibitem[CR12]{CR12}
{A. Chatzistamatiou, K. Rulling},
{\em Hodge-Witt cohomology and Witt-rational singularities,}
{Documenta Math., {\textbf{17}} (2012), 663--781. }



\bibitem[Eke84]{Eke84}
{T. Ekedahl},
{\em On the multiplicative properties of the de Rham-Witt complex. I,}
{Ark. Mat. {\textbf{22}} (1984), no. 2, 185--239.}


\bibitem[FGI$\,^+$05]{FGAex}
{B. Fantechi, L. G{\"o}ttsche, L. Illusie, S. L. Kleiman, N. Nitsure, A. Vistoli},
{\em Fundamental algebraic geometry},
{American Mathematical Society, Providence, RI, 2005}.



\bibitem[Fu11]{Fu11}
{L. Fu},
{\em Etale cohomology theory,}
{Nankai Tracts in Mathematics, {\textbf{13}}. World Scientific Publishing Co. Pte. Ltd., Hackensack, NJ, 2011.}




\bibitem[Jou83]{Jou83}
{J. P. Jouanolou},
{\em Th\'eor\`emes de Bertini et applications,}
{Progress in Mathematics, {\textbf{42}}. Birkhauser Boston, Inc., Boston, MA, 1983.}





\bibitem[Gro64]{Gro64}
{A. Grothendieck},
{\em \'El\'ements de g\'eom\'etrie alg\'ebrique. I. Le langage des sch\'emas,}
{Inst. Hautes \'Etudes Sci. Publ. Math. No. {\textbf{4}} 1960.}




\bibitem[Gro65]{Gro65} 
{A. Grothendieck},
{\em \'El\'ements de g\'eom\'etrie alg\'ebrique. IV. 
\'Etude locale des sch\'emas et des morphismes de sch\'emas. II,}
{Inst. Hautes Etudes Sci. Publ. Math. No. {\textbf{24}} (1965), 231 pp.}




\bibitem[Har67]{Har67}
{R.~Hartshorne},
{\em Local cohomology},
{A seminar given by A. Grothendieck, Harvard University, Fall, 
1961 Lecture Notes in Mathematics, No. {\textbf{41}} Springer-Verlag, Berlin-New York 1967}.


\bibitem[Har77]{Har77}
{R.~Hartshorne},
{\em Algebraic Geometry},
{Grad. Texts in Math., no {\textbf{52}}, Springer-Verlag, NewYork, 1977}.


\bibitem[Haz12]{Haz12}
{M. Hazewinkel},
{\em Formal groups and applications,}
{Corrected reprint of the 1978 original. AMS Chelsea Publishing, Providence, RI, 2012. }



\bibitem[Ill79]{Ill79}
{L.~Illusie}, 
{\em Complexe de de Rham--Witt et cohomologie cristalline},
{Ann. Sci. \'Ecole Norm. Sup., Vol. {\textbf{12}}, Issue 4, (1979), 501--661.}




\bibitem[Kat89]{Kat89}
{K. Kato},
{\em Swan conductors for characters of degree one in the imperfect residue field case. Algebraic K-theory and algebraic number theory (Honolulu, HI, 1987), }
{101--131, Contemp. Math., 83, Amer. Math. Soc., Providence, RI, 1989.}


\bibitem[KMM87]{KMM87}
{Y. Kawamata, K. Matsuda, K. Matsuki},
{\em Introduction to the minimal model program},
{Algebraic geometry, sendai, 1987, pp. 283--360, 
Adv. Stud. Pure Math., vol {\textbf{10}}}.



\bibitem[Kun76]{Kun76}
{E. Kunz},
{\em On Noetherian rings of characteristic $p$},
{Amer. J. Math. {\textbf{98}} (1976), no. 4, 999--1013.}




\bibitem[Mat89]{Mat89}
{H. Matsumura},
{\em Commutative ring theory}, 
{Translated from the Japanese by M. Reid. Second edition. Cambridge Studies in Advanced Mathematics, 
{\textbf{8}}. Cambridge University Press, Cambridge, 1989.}


\bibitem[MR85]{MR85}
{V.~B.~Mehta, A.~Ramanathan}, 
{\em Frobenius splitting and cohomology vanishing for Schubert varieties},
{Ann. of Math. (2), {\textbf{122}} (1985), 27--40}.


\bibitem[Muk13]{Muk13}
{S. Mukai},
{\em Counterexamples to Kodaira's vanishing and Yau's inequality in positive characteristics}, 
{Kyoto J. Math. 53 (2013), no. 2, 515--532.}


\bibitem[Ray78]{Ray78}
{M. Raynaud,} 
{\em Contre-exemple au \lq\lq vanishing theorem'' en caract\'eristique $p>0$}, 
{Ramanujam--a tribute, 273--278,
Tata Inst. Fund. Res. Studies in Math., {\textbf{8}}, Springer, Berlin-New York, 1978. }


\bibitem[Ser79]{Ser79}
{J. P. Serre,} 
{\em Local fields}, 
{Translated from the French by Marvin Jay Greenberg. Graduate Texts in Mathematics, {\textbf{67}}. 
Springer-Verlag, New York-Berlin, 1979.}










\end{thebibliography}
\end{document}